\documentclass[preprints,article,accept,moreauthors]{Definitions/mdpi} 

\firstpage{1} 
\pubvolume{xx}
\issuenum{1}
\articlenumber{5}
\pubyear{2019}
\copyrightyear{2019}

\newcommand{\zoro}[1]{{\textcolor{black}{#1}}}
\newcommand{\yoro}[1]{{\textcolor{black}{#1}}}

\history{Accepted: 17/09/2020; Published: date}

\usepackage{xparse}
\usepackage{chngcntr}
\usepackage{apptools}
\usepackage{amsbsy}
\usepackage{bm}
\usepackage{fixmath}
\usepackage{amssymb}
\usepackage{empheq}

\hypersetup{
colorlinks=true,
urlcolor=blue,
linkcolor=red,
citecolor=green
}

\newcommand{\Fhyp}[5]{\,\mbox{}_{#1}F_{#2}\!\left(
\genfrac{}{}{0pt}{}{#3}{#4};#5\right)}

\newtheorem{thm}[lemma]{Theorem}
\newtheorem{cor}[lemma]{Corollary}
\newtheorem{rem}[lemma]{Remark}

\newtheorem{defn}[lemma]{Definition}
\newtheorem{lem}[lemma]{Lemma}

\makeatletter
\def\eqnarray{\stepcounter{equation}\let\@currentlabel=\theequation
\global\@eqnswtrue
\tabskip\@centering\let\\=\@eqncr
$$\halign to \displaywidth\bgroup\hfil\global\@eqcnt\z@
$\displaystyle\tabskip\z@{##}$&\global\@eqcnt\@ne
\hfil$\displaystyle{{}##{}}$\hfil
&\global\@eqcnt\tw@ $\displaystyle{##}$\hfil
\tabskip\@centering&\llap{##}\tabskip\z@\cr}

\def\endeqnarray{\@@eqncr\egroup
\global\advance\c@equation\m@ne$$\global\@ignoretrue}

\def\@yeqncr{\@ifnextchar [{\@xeqncr}{\@xeqncr[5pt]}}
\makeatother

\newcommand{\N}{{\mathbb N}}

\newcommand{\CC}{{\mathbb C}}
\newcommand{\expe}{{\mathrm e}}

\newcommand{\biQ}{{\mathbold Q}}

\numberwithin{equation}{section}
\numberwithin{corollary}{section}
\numberwithin{remark}{section}
\numberwithin{theorem}{section}
\numberwithin{lemma}{section}
\numberwithin{definition}{section}

\Title{Multi-integral representations for
associated Legendre and Ferrers functions 
}

\Author{Howard S. Cohl$^{1,\dagger}$\orcidA{} and 
Roberto S. Costas-Santos $^{2,\dagger}$\orcidB{}
}

\AuthorNames{Howard S. Cohl, Roberto S. Costas-Santos}

\address{%
$^{1}$ \quad Applied and Computational Mathematics Division,
National Institute of Standards and Technology,
Mission Viejo, CA 92694, USA; howard.cohl@nist.gov\\
$^{2}$ \quad 
Departamento de F\'isica y Matem\'{a}ticas,
Universidad de Alcal\'{a}, c.p. 28871, Alcal\'{a} de Henares, Spain;
rscosa@gmail.com
}

\firstnote{These authors contributed equally to this work.}

\abstract{For the associated Legendre and Ferrers
functions of the first and second kind, we obtain
{new}
multi-derivative and multi-integral 
representation formulas. {The multi-integral
representation formulas that we derive for these
functions generalize some}
{classical
multi-integration formulas.}
As a result of the determination of these formulae,
we compute {some interesting} special values and integral representations
for certain particular combinations of the degree and 
order including the case where there is symmetry and 
antisymmetry for the degree and {order} parameters.
As a consequence of our analysis, we obtain 
{some new} results 
for the associated Legendre function of the second kind {including
parameter values for 
which this 
function is identically zero}.
}


\begin{document}
\section{Introduction}
Using analysis for fundamental solutions of the Laplace equation
on Riemannian manifolds of constant curvature, we have previously 
obtained antiderivatives and integral representations for the
associated Legendre and Ferrers functions of the second kind with 
degree and order equal to within a sign. 
For instance, in 
\cite[Theorem 1]{Cohlhypersphere}, 
using the $d$-dimensional hypersphere with $d=2,3,4,\ldots$, we derive
an antiderivative and an integral representation for the Ferrers 
function of the second kind with order equal to the negative degree.
In 
\cite[Theorem 3.1]{CohlKalII}, 
using the $d$-dimensional hyperboloid model of hyperbolic geometry with 
$d=2,3,4,\ldots$, the authors derived an antiderivative and an 
integral representation for the associated Legendre function of the 
second kind with degree and order equal to each other.
In \cite{Cohlhypersphere,CohlKalII}, the antiderivatives
and integral representations were restricted to values of the degree 
and order $\nu$ such that $2\nu$ is an integer.
\yoro{One of the goals of} this paper is to generalize \yoro{some integral representation} results presented in
\cite{Cohlhypersphere,CohlKalII} for the associated Legendre and 
Ferrers functions of the first and second kind, and to extend them 
such that the degree and order are no longer subject to the above 
restriction. Our integral representations are 
consistent with known special values for the associated Legendre 
and Ferrers functions of the first kind when the order is equal 
to the negative degree.

\zoro{
\phantom{XXX}The multi-integrals presented in this paper 
(Theorems \ref{critthmOlverQ},
\ref{assLegmultint},
\ref{assLegmultint-2},
\ref{repintFerPzx},
\ref{repintFerPzxb} and 
\ref{repintFerPzxb})
generalize multi-integrals 
for  arbitrary order
$(\mu)$,
which have
appeared previously in the literature
(see for instance,
\cite[Section 14.6(ii)]{NIST:DLMF},
\cite[(8.14.17-18)]{Abra} and the earliest
appearance we have found
\cite[p.~149]{Erdelyi}).
The multi-integrals for the 
Ferrers function of the second 
kind 
(Theorems \ref{CorFerQnummu},
\ref{CorFerQnumu},
\ref{repintFerQzx}
and
\ref{repintFerQzxb})
also produce generalized results for
arbitrary order $(\mu)$.
In fact, the specialization
of these multi-integrals
for $\mu=0$ has not 
appeared in the literature.
}

\zoro{
\phantom{XXX}Applications of the work contained in this manuscript 
include any of the many different areas in which Legendre and Ferrers functions arise, which include a very large number of disciplines. The associated
Legendre and Ferrers functions treated in this paper, 
e.g., $Q_\nu^\mu(\cosh r)/\sinh^\mu r$, appear as fundamental
solutions of the Laplace and Helmholtz equation on Riemannian manifolds
of constant curvature (see e.g., \cite[Section 3.3]{CohlDangDunster}). Associate Legendre and Ferrers functions appear in any place where harmonic analysis needs 
to be performed on the surface of a sphere or on an oblate or prolate spheroid. These are analogs of the $1/r$
potential in Euclidean space for Riemannian spaces of constant curvature.
Therefore results such as we derive below will be important
when studying global analysis of these fundamental 
solutions for higher powers of the Laplacian or Helmholtz operators. These multi-integration results provide an algorithm for computing fundamental solutions of much larger powers
of the Laplace-Beltrami operator on these spaces. A survey of applications of associated Legendre and Ferrers functions is given in 
\cite[Sections 14.30-31]{NIST:DLMF}.
This points to harmonic analysis
on the surface of spheres, oblate and prolate spheroids, and
circular toroids. Other applications
include the Mehler-Fock transforms, high frequency atomic and molecular scattering, quantum direct and exchange
Coulomb interaction, Newtonian gravity, etc.
Also, the special cases treated in this paper, many
of which are also not well-known, provide for a beautiful illumination of this classical subject.
}

\zoro{
The critical aspects that allows for
proofs of the results
contained in this paper is that fact that we are able to obtain differentiation/integration
properties of associated Legendre and Ferrers functions such that
the order $(\mu)$ of these functions are either raised or
lowered by
integral amounts. Furthermore, the degree $(\nu)$ of these functions
are not at all affected by the differentiation/integration. In fact, these
differentiation/integration properties (Remarks \ref{rem39},
\ref{rem315}, \ref{renderFerP},
\ref{rem47},
\ref{rem415},
\ref{rem419},
\ref{rem424} and
\ref{rem428})
are not well-known in the literature \yoro{for} these functions.
}

{Note that in the process of our derivations,}
we also obtained some nice results for the associated
Legendre function of the second kind with degree ${\nu=-\frac32-n}$.
{This included the} full Gauss hypergeometric dependence
{and} large argument 
asymptotics {which is 
crucial for establishing 
Theorem \ref{critthmOlverQ}}. 
We are also
able to find their zeros which occur when 
$\mu=\pm(\frac12+k)$, $n,k\in\mathbb N_0$, $n\ge k$ {(see  Corollary \ref{zerosLegQ})}.

\medskip
\section{{Preliminaries}}

Throughout this paper we adopt the following set notations:
$\mathbb N_0:=\{0\}\cup\mathbb N=\{0, 1, 2, 3, \ldots\}$, and
we use the set $\mathbb C$ which represents the complex
numbers.
As is the common convention for associated Legendre 
functions \cite[(8.1.1)]{Abra}, for any expression of the 
form $(z^2-1)^\alpha$, read this as 
$(z^2-1)^\alpha:=(z+1)^\alpha(z-1)^\alpha$,
for any fixed $\alpha\in\mathbb C$ and $z\in\mathbb C\setminus(-\infty,1]$.
In this paper, we will use
the Gauss hypergeometric 
function ${}_2F_1$ which 
can be defined in terms of the following infinite series
as \cite[(2.1.5)]{AAR} 
\[
\Fhyp21{a,b}{c}{z}:=
\sum_{n=0}^\infty \frac{(a)_n(b)_n}{(c)_n}\frac{z^n}{n!},
\]
where $c\not\in-\mathbb N_0$,
and elsewhere on $z\in\mathbb C\setminus(1,\infty)$ by analytic
continuation; {where the Pochhammer symbol (rising factorial)
is defined by
\begin{equation}
(z)_n:=\prod_{i=1}^n(z+i-1),
\label{Pochdef}
\end{equation}
}
with $n\in\mathbb N_0$. 
Note that
for all $n\in\mathbb N_0$, $z\not\in-\mathbb N_0$,
one has
\begin{equation}
(z)_n=\frac{\Gamma(z+n)}{\Gamma(z)}, \quad
\Gamma(z-n)=\frac{(-1)^n\Gamma(z)}
{(-z+1)_n}.
\label{Gamrefneg}
\end{equation}
We will also need the binomial theorem \cite[(15.4.6)]{NIST:DLMF}
\begin{equation}
\Fhyp21{a,b}{b}{z}=(1-z)^{-a},
\label{binom}
\end{equation}
{
Euler's transformation \eqref{Euler},
\cite[(15.8.1)]{NIST:DLMF}
\begin{equation}\label{Euler}
\Fhyp21{a,b}{c}{z}=(1-z)^{c-a-b}\Fhyp21{c-a,c-b}{c}{z},
\end{equation}
and the Gauss sum \cite[(15.4.20)]{NIST:DLMF}
\begin{equation}
\Fhyp21{a,b}{c}{1}=
\frac{\Gamma(c)\Gamma(c-a-b)}
{\Gamma(c-a)\Gamma(c-b)},
\label{Gauss}
\end{equation}
for $\Re(c-a-b)>0$.
}
We will also use the {generalized} hypergeometric
function
\[
\Fhyp32{a,b,c}{d,e}{z}
:=
\sum_{n=0}^\infty \frac{(a)_n(b)_n(c)_n}{(d)_n(e)_n}\frac{z^n}{n!},
\]
where $d,e\not\in-\mathbb N_0$.
\medskip
\noindent 
We now produce a lemma for the Gauss hypergeometric
function which will be useful 
in our analysis of antiderivatives for associated
Legendre functions of the second kind below.
\begin{lem}
\label{derivlem}
Let $z\in\mathbb C\setminus[0,\infty)$.
Then
\begin{equation}
\frac{d}{dz}
\frac{1}{z^{\nu+\mu+1}}
\Fhyp21{\frac{\nu+\mu+1}{2},\frac{\nu+\mu+2}{2}}{\nu+\frac32}{\frac{1}{z^2}}
=\frac{-(\nu+\mu+1)}{z^{\nu+\mu+2}}
\Fhyp21{\frac{\nu+\mu+2}{2},\frac{\nu+\mu+3}{2}}{\nu+\frac32}{\frac{1}{z^2}}.
\label{derivhyp21form}
\end{equation}
\end{lem}
\begin{proof}
Differentiating the left-hand side
of \eqref{derivhyp21form} using
the chain rule and \cite[(15.5.1)]{NIST:DLMF}
\begin{equation}
\frac{d}{dz}\Fhyp21{a,b}{c}{z}
=\frac{ab}{c}
\Fhyp21{a+1,b+1}{c+1}{z},
\label{derhyp}
\end{equation}
one produces an expression involving the sum of 
two Gauss hypergeometric functions. One can then 
use the following Gauss relations for contiguous
hypergeometric functions {\cite[p. 58]{Erdelyi}}
\begin{equation}
z\,{}_2F_1\left(\begin{array}{c}a+1,b+1\\[2mm]
c+1\end{array};z\right)=\frac{c}{a-b}\left[{}_2F_1
\left(\begin{array}{c}a,b+1\\
c\end{array};z\right)-{}_2F_1\left(\begin{array}{c}a+1,b\\
c\end{array};z\right)\right],
\label{sechypcontig}
\end{equation}
and \cite[(15.5.12)]{NIST:DLMF}
\begin{equation}
{}_2F_1\left(\begin{array}{c}a,b+1\\[0.2cm]
c\end{array};z\right)=\frac{b-a}{b}\,{}_2F_1\left(\begin{array}{c}a,b
\\[0.2cm]c\end{array};z\right)
+\frac{a}{b}\,{}_2F_1\left(\begin{array}{c}a+1,b\\[0.2cm]c\end{array};z\right),
\label{thirdhypcontig}
\end{equation}
to obtain the following formula
\begin{equation}
\frac{d}{dz}\frac{1}{z^{a+b-\frac12}}
\Fhyp{2}{1}{a,b}{c}{\frac{1}{z^2}}=
\frac{a-b+\frac12}{z^{a+b+\frac12}}
\Fhyp{2}{1}{a,b}{c}{\frac{1}{z^2}}
-\frac{2a}{z^{a+b+\frac12}}
\Fhyp{2}{1}{a+1,b}{c}{\frac{1}{z^2}}.
\label{derivQhyps}
\end{equation}
Since for the Gauss hypergeometric 
function on the left-hand side of 
\eqref{derivhyp21form}, $a-b+\frac12=0$, 
so the first term on the right-hand side
of \eqref{derivQhyps}
vanishes and the lemma follows.
\end{proof}

{
\begin{defn}
Let $z\in\CC$, $a,b\in\CC\cup\{-\infty,\infty\}$.
Define the following notations for $n$th iterated integrals 
of the functions 
$f(z;{\bf a})$, $g(z;{\bf b})$, 
respectively, 
\begin{eqnarray}
&&\hspace{-0.5cm} 
\int_z^b \cdots\int_z^b f(w;{\bf a}) (dw)^n 
:=\int_z^b \left[\int_{w_{n-1}}^b\cdots\,
\left[\int_{w_2}^b \left[\int_{w_1}^b f(w;{\bf a})\, dw\,\right]
dw_1\,\right]\cdots\,dw_{n-2}\,\right]dw_{n-1},\\
&&\hspace{-0.5cm} 
\int_a^z\cdots\int_a^z g(w;{\bf b}) (dw)^n 
:=\int_a^z\left[\int_a^{w_{n-1}}\cdots\,
\left[\int_a^{w_2}\left[\int_a^{w_1}g(w;{\bf b})\,dw\,\right]
dw_1\,\right]\cdots\,dw_{n-2}\,\right]dw_{n-1},
\end{eqnarray}
where $w_0:=w$, $w_{n}:=z$, and 
${\bf a}$, ${\bf b}$,
are sets of fixed parameters.
\end{defn}
}

Another useful result we are going to use 
often along this work is the following.
\begin{lem}
Let $n\in\N_0$, \yoro{$a, x,\mu\in \mathbb C$}, and let $f^\mu$ be a 
function such that
\begin{equation} \label{eq:fmud}
\frac {d}{dz} f^\mu(z)=\lambda_\mu f^{\mu\pm 1}(z),
\end{equation}
{where $\lambda_\mu\in \mathbb \CC^\ast$. 
\yoro{Then,} the following identity holds:}
\[
\int_a^x\cdots \int_a^x f^\mu(w)
(dw)^n
=
\frac 1{\lambda_{\mu\mp 1}\cdots \lambda_{\mu\mp n}}
\sum_{k=n}^\infty \frac{\lambda_{\mu\mp n}\cdots 
\lambda_{\mu\mp n\pm(k-1)} f^{\mu\mp n\pm k}(a)(x-a)^k}{k!}.
\]
\label{Robertolem}
\end{lem}

\begin{proof}
We are going to prove the result by induction on $n$.
\yoro{The $n=0$ case is direct taking into account the Taylor expansion of 
$f$ at $x=a$.}
If $n=1$ then
\[
\int_a^x f^\mu(w)dw=\frac 1{\lambda_{\mu\mp1}} 
\int_a^x 
\left(\frac {d}{dw} f^{\mu\mp 1}(w)\right)\, dw=
\frac 1{\lambda_{\mu\mp1}}\left(f^{\mu\mp 1}(x)-
 f^{\mu\mp 1}(a)\right).
\]
By using the Taylor expansion of $f^{\mu\mp 1}(x)$ at $x=a$
and \eqref{eq:fmud}, the result follows for the $n=1$ case.
Assuming the result holds for $n$, let us prove the 
identity for the $n+1$ case:
\begin{eqnarray*}
\int_a^x\cdots \int_a^x f^\mu(w)
(dw)^n
&=&
\frac 1{\lambda_{\mu\mp 1}}\int_a^x\cdots \int_a^x 
\left(f^{\mu\mp 1}(w_1)-f^{\mu\mp 1}(a)\right) 
dw_1\cdots d w_{n}
\\
&=&
\frac 1{\lambda_{\mu\mp 1}\cdots\lambda_{\mu\mp (1+n)}}
\sum_{k=n+1}^\infty \frac{\lambda_{\mu\mp (n+1)}\cdots 
\lambda_{\mu\mp (n+1)\pm (k-1)} f^{\mu\mp (n+1)\pm k}(a)(x-a)^k}{k!},
\end{eqnarray*}
where we have used induction and the basic properties of integrals.
Hence the result follows. 
\end{proof}

\section{Associated Legendre functions of the first and second kind}
Associated Legendre functions 
(and Ferrers functions) are 
those Gauss hypergeometric functions
which satisfy a quadratic transformation
(see \cite[Sections 15.8(iii-iv)]{NIST:DLMF}).
In the following sections we
will derive derivative, antiderivative,
and integral representations for associated
Legendre (and Ferrers) functions of the
first and second \yoro{kinds}which to the best
of our knowledge have not appeared in
the classical literature of these
highly applicable special functions
of applied and pure mathematics.

\medskip

The associated Legendre function of the first kind
 $P_\nu^\mu :\mathbb C\setminus(-\infty,1]\to\mathbb C$ 
is defined as \cite[(14.3.6)]{NIST:DLMF}
\begin{equation}
P_\nu^\mu(z)=
\frac{1}{\Gamma(1-\mu)}
\left(
\frac{z+1}{z-1}
\right)^{\mu/2}
\Fhyp21{-\nu,\nu+1}{1-\mu}{\frac{1-z}{2}}.
\label{defLegP}
\end{equation}
Starting with 
\eqref{defLegP}, setting \yoro{$\mu\mapsto -\mu$,} and applying
\eqref{Euler},
\yoro{another useful hypergeometric
representation} for the associated Legendre function of
the first kind \yoro{can be obtained}, namely
\begin{equation}
P_\nu^{-\mu}(z)=
\frac{(z^2-1)^{\frac{\mu}{2}}}
{2^\mu\Gamma(\mu+1)}
\Fhyp21{\nu+\mu+1,-\nu+\mu}{1+\mu}
{\frac{1-z}{2}}.
\label{EulerLegP}
\end{equation}
The associated Legendre function of the second 
kind $\biQ_\nu^\mu:\mathbb C\setminus(-\infty,1]\to\mathbb C$
can be defined in terms of the Gauss hypergeometric
function
as \cite[
(14.3.10) and Section 14.21]{NIST:DLMF}
\begin{equation}
\biQ_\nu^\mu(z):=\frac{\sqrt{\pi}
(z^2-1)^{\mu/2}}
{2^{\nu+1}\Gamma(\nu+\frac32)z^{\nu+\mu+1}}
\Fhyp21{\frac{\nu+\mu+1}{2},\frac{\nu+\mu+2}{2}}{
\nu+\frac32}{\frac{1}{z^2}},
\label{Qdefn}
\end{equation}
for $|z|>1$ and, by analytic continuation of the Gauss hypergeometric 
function, elsewhere on $z\in\mathbb C\setminus(-\infty,1]$.
\begin{rem}
\label{Hobsonnot}
The normalized notation $\biQ_\nu^\mu(z)$
is due to Olver \cite[p.~178]{Olver}
and is defined in terms of the more
commonly appearing Hobson notation for
the associated Legendre function of the
second kind $Q_\nu^\mu(z)$ as follows
\cite[(14.3.10)]{NIST:DLMF}
\begin{equation}
\biQ_\nu^\mu(z)=\frac{\expe^{-i\pi\mu}}{\Gamma(\nu+\mu+1)}Q_\nu^\mu(z).
\end{equation}
See \cite[Section 14.1]{NIST:DLMF} for more
information on commonly appearing
notations for the associated Legendre and
Ferrers functions.
\end{rem}

\begin{rem}
\yoro{Note that the following
algebraic special cases for the associated
Legendre function of the second}\\ 
\yoro{kind, hold for $\mu=\nu+1,\nu+2$, namely}
\begin{equation}
\biQ_\nu^{\nu+1}(z)
=\frac{\sqrt{z^2-1}}{z}
\biQ_\nu^{\nu+2}(z)
=\dfrac{\sqrt \pi}{2^{\nu+1}\Gamma(\nu+\frac{3}{2})(z^2-1)^{\frac{\nu+1}2}},
\label{Specialmunu12}
\end{equation}
where we have used \eqref{Qdefn}
and \eqref{binom}.
Furthermore algebraic expressions
for $\biQ_\nu^{\nu+n}$ for all
$n\in\N$ are obtainable from the order
recurrence relation for associated 
Legendre functions of the second kind
(cf.~\cite[(14.10.6)]{NIST:DLMF})
\begin{equation}
\biQ_\nu^{\nu+n+2}(z)=
\frac{2(\nu+n+1)z}{(2\nu+n+2)\sqrt{z^2-1}}
\biQ_\nu^{\nu+n+1}(z)-
\frac{n}{2\nu+n+2}
\biQ_\nu^{\nu+n}(z).
\end{equation}
Since $\biQ_\nu^{\nu+2}(z)$ is proportional
to $\biQ_\nu^{\nu+1}(z)$, then all 
$\biQ_\nu^{\nu+m}(z)$ is proportional
to $\biQ_\nu^{\nu+1}(z)$ for all $m\ge 2$.
\end{rem}

We now present some theorems related to the behavior
of the associated Legendre function of the second
kind with degree $\nu=-\frac32-n\in\{-\frac32,-\frac52,\ldots\}$, 
$n\in\mathbb N_0$ and its corresponding asymptotics
as $z\to\infty$. This will be useful in our further
analysis below.

\begin{thm}
Let $z\in\mathbb C\setminus(-\infty,1]$,
$\mu\in\mathbb C$, $\nu=-\frac32-n\in\{-\frac32,-\frac52,\ldots\}$, $n\in\mathbb N_0$. Then
\begin{equation}
\biQ_{-\frac32-n}^\mu(z)=
\frac{(-1)^n\sqrt{\pi} (\mu^2-\frac14)(\frac32+\mu)_n(\frac32-\mu)_n(z^2-1)^{\frac{\mu}{2}}}
{2^{n+\frac32}(n+1)!\,z^{n+\frac32+\mu}}
\Fhyp21{
\frac34+\frac{\mu+n}{2},
\frac54+\frac{\mu+n}{2}}{n+2}{\frac{1}{z^2}}.
\label{assLegQnum32}
\end{equation}
\end{thm}
\begin{proof}
Start with \eqref{Qdefn} then let $\nu=-\frac32-n$,
followed by the application of 
\cite[First equation in Section 15.2(ii)]{NIST:DLMF},
\[
\lim_{c\to-n}\frac{1}{\Gamma(c)}
\Fhyp21{a,b}{c}{z}=
\frac{(a)_{n+1}(b)_{n+1}}
{(n+1)!}z^{n+1}
\Fhyp21{a+n+1,b+n+1}{n+2}{z},
\]
and the Pochhammer symbol identity for $\nu\in\mathbb C$, $n\in\mathbb N_0$, 
\[
\left(\nu-\frac{n}{2}\right)_n
\left(\nu+\frac12-\frac{n}{2}\right)_n=
\frac{(-1)^n}{2^{2n}}(2\nu)_n(-2\nu+1)_n,
\]
which follows from the duplication theorem for 
gamma functions then \eqref{Gamrefneg}.
\end{proof}

The following corollary is an interesting
side-effect of the above theorem which produces 
zeros for the associated Legendre function
of the second kind.
\begin{cor}
Let $z\in\mathbb C\setminus(-\infty,1]$, $\nu=-\frac32-n\in\{-\frac32,-\frac52,\ldots\}$, $\mu=\pm(\frac12+k)$, $n,k\in\mathbb N_0$. Then
\begin{equation}
\biQ_\nu^{\mu}(z)=\biQ_{-\frac32-n}^{\pm(\frac12+k)}(z)=0,
\end{equation}
for all $n\ge k$.
\label{zerosLegQ}
\end{cor}
\begin{proof}
Simple examination of the factor multiplying the
Gauss hypergeometric function
in \eqref{assLegQnum32} produces the result.
\end{proof}

\begin{rem}
Note that the zeros for associated Legendre
functions for the $k=0$ case in 
Corollary \ref{zerosLegQ} 
is clear from the special value
\cite[(14.5.17)]{NIST:DLMF}
\[
\biQ_\nu^{\pm\frac12}(\cosh\xi)
=\sqrt{\frac{\pi}{2\sinh\xi}}
\frac{\exp(-(\nu+\frac12)\xi)}{\Gamma(\nu+\frac32)}.
\]
\end{rem}

We now give a result which produces the large argument
asymptotics for the associated Legendre function
of the second kind when the degree $\nu=-\frac32-n$,
$n\in\mathbb N_0$.
\begin{lem}
\label{newaLegQasymp}
Let $z\in\mathbb C\setminus(-\infty,1]$,
$\mu\in\mathbb C$, $\nu=-\frac32-n\in\{-\frac32,-\frac52,\ldots\}$, $n\in\mathbb N_0$. Then
\begin{equation}
\biQ_\nu^\mu(z)\sim
\frac{(-1)^{-\nu-\frac32}\sqrt{\pi}\,2^\nu (\mu^2-\frac14)
\Gamma(\mu-\nu)\Gamma(-\mu-\nu)z^\nu}
{\Gamma(\frac12-\nu)\Gamma(\frac32+\mu)\Gamma(\frac32-\mu)}.
\end{equation}
Equivalently, 
\begin{equation}
\biQ_{-\frac32-n}^\mu(z)
\sim
\frac{(-1)^n\sqrt{\pi}\,(\mu^2-\frac14)(\frac32+\mu)_n
(\frac32-\mu)_n}
{2^{n+\frac32}(n+1)!z^{n+\frac32}}.
\end{equation}
\end{lem}
\begin{proof}
The result follows by starting with 
\eqref{assLegQnum32} and examining
{its leading term} behavior as $z\to\infty$.
\end{proof}
\subsection{The associated Legendre function of the second kind}
We now compute some antiderivatives and 
integral representations 
for associated Legendre functions
of the second kind. This also includes
some nice limits and specializations.


\medskip
\zoro{
\begin{rem}
Note the following expression can be obtained by 
using the definition \eqref{Qdefn}
and Lemma \ref{derivlem}
for $z\in\mathbb C\setminus(-\infty,1]$, $\nu,\mu\in\mathbb C$:
\begin{equation}
\dfrac{d}{dz}\dfrac {\biQ_\nu^\mu(z)}{(z^2-1)^{\frac \mu 2}}=
\dfrac{-(\nu+\mu+1)}{(z^2-1)^{\frac {\mu+1}2}}\biQ_{\nu}^{\mu+1}(z).
\label{dervboldQ}
\end{equation}
From this formula 
the following antiderivative is obtained:
\begin{equation}
\int\frac{\biQ_\nu^{\mu}(z)}
{(z^2-1)^{\frac{\mu}{2}}}
dz
=
\frac{-\biQ_\nu^{\mu-1}(z)}
{(\nu+\mu)(z^2-1)^{\frac{\mu-1}{2}}}+C,
\label{antiderbodLegQ}
\end{equation}
where $C$ is an arbitrary constant.
\end{rem}
}

\begin{thm}
\label{thmintonceLegQ}
Let $z\in\mathbb C\setminus(-\infty,1]$, $\nu,\mu\in\mathbb C$,
such that $\Re(\nu+\mu+1)>0$.
Then
\begin{equation}
\biQ_\nu^{\mu}(z)=(\nu+\mu+1)(z^2-1)^\frac{\mu}{2}
\int_z^\infty\frac{\biQ_\nu^{\mu+1}(w)}{(w^2-1)^{\frac{\mu+1}{2}}}dw.
\label{singleintboldLegQ}
\end{equation}
\end{thm}
\begin{proof}
Taking the limit of the antiderivative 
\eqref{antiderbodLegQ}
evaluated at the endpoints of integration 
using the large argument asymptotics 
\cite[(14.8.15)]{NIST:DLMF}
\[
\biQ_\nu^\mu(z)\sim\frac{\sqrt{\pi}}{\Gamma(\nu+\frac32)(2z)^{\nu+1}},\quad \nu\not\in\{-\tfrac32,-\tfrac52,-\tfrac72,\ldots\},
\]
and Lemma \ref{newaLegQasymp} which shows that 
$\biQ_\nu^\mu(z)\to0$ as $z\to\infty$ for
$\nu\in\{-\tfrac32,-\tfrac52,-\tfrac72,\ldots\}$
as well. Therefore the integral is convergent as indicated
which completes the proof.
\end{proof}

\begin{rem}
\label{rem39}
Iterating \eqref{dervboldQ}, then using induction
with \eqref{Pochdef}, 
\yoro{the following order-shift derivative formula for}
\yoro{the associated Legendre 
function of the second kind, namely for
$z\in\mathbb C\setminus(-\infty,1]$,
$n\in\mathbb N_0$, $\nu,\mu\in\mathbb C$, holds:}
\begin{equation}
\dfrac{d^n}{dz^n}\dfrac{\biQ_\nu^\mu(z)}
{(z^2-1)^{\frac \mu 2}}=
\dfrac{(-1)^n (\nu+\mu+1)_n}{(z^2-1)^{\frac {\mu+n}2}}\biQ_{\nu}^{\mu+n}(z).
\label{oneintLegQ}
\end{equation}
\end{rem}

\begin{thm}
\label{critthmOlverQ}
Let $n\in\N_0$, $z\in\mathbb C\setminus(-\infty,1]$, $n\in\mathbb N_0$, $\nu,\mu\in\mathbb C$, 
such that $\Re(\nu+\mu-n+1)>0$. Then
\begin{equation}
\int_z^\infty \cdots\int_z^\infty 
\dfrac {\biQ_\nu^\mu(w)}{(w^2-1)^{\frac \mu2}} (dw)^n=
\frac{(-1)^n\biQ_\nu^{\mu-n}(z)}{ (-\nu-\mu)_n(z^2-1)^{\frac{\mu-n}2}}
.
\label{intrepzeromQ}
\end{equation}
\end{thm}
\begin{proof}
Iterating Theorem 
\ref{thmintonceLegQ} 
with \eqref{Pochdef}
using induction 
with \eqref{Pochdef}
completes the proof.
\end{proof}

\zoro{
\begin{rem}
It is clear that Theorem \ref{critthmOlverQ}
is a generalization of cf.~\cite[(14.6.8)]{NIST:DLMF}
\begin{equation}
\biQ_\nu^{-n}(z)=(-1)^n (-\nu)_n(z^2-1)^{-\frac{n}2}
\int_z^\infty \cdots\int_{z}^\infty 
\biQ_\nu(w)(dw)^n,
\end{equation}
by considering the $\mu=0$ specialization in \eqref{intrepzeromQ},
using \eqref{Gamrefneg} 
and Hobson's notation
(see Remark \ref{Hobsonnot}).
\end{rem}
}

\begin{rem}
An antiderivative
of an algebraic function (essentially in terms of reciprocal powers 
of the hyperbolic sine function)
expressed as the associated Legendre 
function of the second kind with order and degree equal to each other
\yoro{can be obtained}. 
This is accomplished
by starting with \eqref{antiderbodLegQ}
and setting $\mu=\nu+1$, then using
\eqref{Specialmunu12}. \yoro{This 
produces the specialized
antiderivative,} namely for
\label{antider}
$z\in\mathbb C\setminus(-\infty,1],$ 
$\nu\in\mathbb C\setminus\{-\frac12,-\frac32,-\frac52,\ldots\}$,
\begin{eqnarray}
\int\frac{dz}{(z^2-1)^{\nu+1}}&=&\frac{-1}{(2\nu+1)z^{2\nu+1}}
\,{}_2F_1\left(\begin{array}{c}\nu+\frac12,\nu+1\\[0.2cm]
\nu+\frac32\end{array};\frac{1}{z^2}\right)+C\nonumber
=
\frac{-2^{\nu}\Gamma(\nu+\frac12)}
{\sqrt{\pi}(z^2-1)^{\frac{\nu}{2}}}\biQ_\nu^\nu(z)+C,
\label{antiderivativecosh2r}
\end{eqnarray}
where $C$ is an arbitrary constant.
\label{thmantifirst}
\end{rem}
\medskip

A straightforward consequence of the antiderivative
(\ref{antiderivativecosh2r}) is the following
integral representation for the associated 
Legendre function of the 
second kind with degree and order equal to 
each other.

\begin{cor}
Let $\nu\in\mathbb C$ such that $\Re\nu>-\frac12$, $z\in\mathbb C\setminus(-\infty,1]$.
Then
\begin{equation}
\biQ_\nu^\nu(z)=
\biQ_\nu^{-\nu}(z)=
\frac
{\sqrt{\pi}
(z^2-1)^{\frac{\nu}{2}}}{2^\nu\Gamma(\nu+\frac12)}
\int_{z}^\infty \frac{dw}{(w^2-1)^{\nu+1}}.
\label{Qnunuz}
\end{equation}
\end{cor}
\begin{proof}
Evaluating the 
antiderivative
Theorem \ref{antider}
at the endpoints of integration
and taking advantage
of \cite[(14.9.14)]{NIST:DLMF}
\[
\biQ_\nu^{-\mu}(z)
=\biQ_\nu^{\mu}(z),
\]
completes the proof.
\end{proof}


\subsection{The associated Legendre function of the first kind}
An integral
representation for the associated
Legendre function of the first 
kind by applying the Whipple
formulae to \eqref{oneintLegQ} \yoro{can be obtained}.
However, this integral representation
shifts the degree of the associated
Legendre function of the first kind
$\nu$ by unity
instead of shifting the order by
unity.

\begin{cor}
Let $z\in\mathbb C\setminus(-\infty,1]$,
$\nu,\mu\in\mathbb C$.
Then
\begin{equation}
P_{\nu-1}^{-\mu}(z)=(\nu+\mu)
(z^2-1)^{-\frac{\nu}{2}}
\int_1^z 
\frac{P_\nu^{-\mu}(w)}
{(w^2-1)^{\frac{\nu+2}{2}}}
dw.
\end{equation}
\end{cor}
\begin{proof}
Apply the Whipple formula 
\cite[(14.9.16)]{NIST:DLMF}
\begin{equation}
\biQ_\nu^\mu(z)=\sqrt{\frac{\pi}{2}}
(z^2-1)^{-1/4}
P_{-\mu-1/2}^{-\nu-1/2}\left(\frac{z}{\sqrt{z^2-1}}\right),
\label{Whipple}
\end{equation}
to the associated Legendre functions
of the second kind on the
left and right-hand sides of
\eqref{oneintLegQ}
followed by the application of the 
involution
\cite[Section 2]{Cohlonpar}
$\zeta(z):=\log\coth\frac{z}{2}$ and making a change of variables
$w=\zeta/\sqrt{\zeta^2-1}$ completes the 
proof.
\end{proof}

The following integral representation
\yoro{can be derived}
by applying Whipple's formulae to our integral 
representation for the associated Legendre function of 
the second kind. We are able to obtain 
an integral representation for the associated Legendre function of 
the first kind which shifts the order by an integer value, 
similar to \eqref{intrepzeromQ}. 
This is achieved by deriving a corresponding derivative formula as follows.

\begin{rem}
\label{rem315}
If you divide both sides of \eqref{EulerLegP}
by $(z^2-1)^{\frac{\mu}{2}}$ and differentiate with respect
to $z$, \yoro{by} using \eqref{derhyp}, 
the unit increments
of the parameters of the Gauss hypergeometric function can 
be absorbed in the order $(\mu)$ of the associated
Legendre function of the first kind.
Let $n\in\mathbb N_0$, $z\in\mathbb C\setminus(-\infty,1]$, 
$\nu,\mu\in\mathbb C$.
Then
\begin{equation}
\frac{d^n}{dz^n}
\frac{P_\nu^{-\mu}(z)}{(z^2-1)^{\frac{\mu}{2}}}
=\frac{(-1)^n (\nu+\mu+1)_n(\mu-\nu)_n}
{(z^2-1)^{\frac{\mu+n}{2}}}
P_\nu^{-\mu-n}(z).
\label{derivnPnummu}
\end{equation}
\end{rem}

\noindent From the above result
integral representations \yoro{can be obtained}
through repeated integration.
For instance, the single integral result
is given as follows.

\begin{cor}
Let $z\in\mathbb C\setminus(-\infty,1]$,
$\nu,\mu\in\mathbb C$.
Then
\begin{equation}
\int_1^z\frac{P_\nu^{-\mu}(w)}{(w^2-1)^{\frac{\mu}{2}}}dw
=\frac{1}{(\nu+\mu)(\nu-\mu+1)}
\left(\frac{P_\nu^{-\mu+1}(z)}{(z^2-1)^{\frac{\mu-1}{2}}}
-\frac{1}{2^{\mu-1}\Gamma(\mu)}\right).
\end{equation}
\end{cor}
\begin{proof}
In order to derive this result,
\yoro{after applying} the fundamental theorem
of calculus for some continuous function $f$
on $[a,b]$,
\begin{equation}
\int_a^b f\,'(t)dt=f(b)-f(a),
\label{funthmcalc}
\end{equation}
and taking advantage of
\cite[(14.8.7)]{NIST:DLMF}
\begin{equation}
\label{limLegPone}
\lim_{z\to1^{+}}
\frac
{P_\nu^{-\mu}(z)}
{(z^2-1)^{\frac{\mu}{2}}}
=
\frac{1}{2^{\mu}\Gamma(\mu+1)},
\end{equation}
\yoro{this} completes the proof.
\end{proof}
\noindent The above result \yoro{can be generalized}
by repeatedly integrating the above
formula. 

\begin{thm}
\label{assLegmultint}
Let $n\in\mathbb N_0$, $z\in\mathbb C\setminus(-\infty,1]$,
$\nu,\mu\in\mathbb C$.
Then
\begin{eqnarray}
\hspace{-10mm}\int_1^z\cdots\int_1^z\frac{P_\nu^{-\mu}(w)}
{(w^2-1)^{\frac{\mu}{2}}}(dw)^n
&=&\frac{1}{(-\nu-\mu)_n(\nu-\mu+1)_n}\Biggl(
\frac{(-1)^n P_\nu^{-\mu+n}(z)}
{ (z^2-1)^{\frac{\mu-n}{2}}}\nonumber\\
&&-\frac
{(-\mu)_n}
{2^{\mu-n}\Gamma(\mu+1)}
\sum_{k=0}^{n-1}
\frac{(\nu+\mu+1-n)_k(\mu-\nu-n)_k}
{k!(\mu+1-n)_k}
\left(\frac{1-z}{2}\right)^k\Biggr)
\label{twosumLegP}\\
&=&\frac{(z-1)^n}{2^\mu n!\Gamma(\mu+1)}
\Fhyp32{\nu+\mu+1,\mu-\nu,1}{\mu+1,n+1}{\frac{1-z}{2}}.
\label{assLegtwosum}
\end{eqnarray}
\label{multintLegP}
\end{thm}

\begin{proof}
Repeated integration of 
\eqref{derivnPnummu} while
noting 
\eqref{limLegPone},
and using induction 
with \eqref{Pochdef}
derives the two sum 
expression \eqref{twosumFerP}.
By rewriting the associated Legendre function of the first
kind on the right-hand side of \eqref{assLegtwosum} in terms of the Gauss 
hypergeometric representation
\eqref{EulerLegP}, 
the finite sum term cancels the first $n$ 
terms of the $k$ sum, and rewriting 
the resulting expression
shows it can be written 
in terms of a nonterminating ${}_3F_2$.
This completes the proof.
\end{proof}

\zoro{
\begin{rem}
It is clear that Theorem \ref{assLegmultint} is a generalization of  \cite[(14.6.7)]{NIST:DLMF}
\begin{equation}
P_\nu^{-n}(z)=(z^2-1)^{-\frac{n}{2}} 
\int_1^z\cdots\int_1^z P_\nu(w)(dw)^n,
\end{equation}
by considering the specialization 
$\mu=0$ in \eqref{twosumLegP}
which follows by using \eqref{Gamrefneg} and \cite[(14.9.13)]{NIST:DLMF}
\begin{equation}
P_\nu^{-n}(z)=\frac{\Gamma(\nu-n+1)}{\Gamma(\nu+n+1)}P_\nu^n(z),\quad n\in\mathbb N_0.
\end{equation}
\end{rem}
}

\begin{rem}
Note the special value
\cite[(14.5.19)]{NIST:DLMF}
\begin{equation}
P_\nu^{-\nu}(z)=\frac{(z^2-1)^{\frac{\nu}{2}}}{2^\nu\Gamma(\nu+1)}.
\label{legPnumnu}
\end{equation}
\end{rem}

Using this special value and connection properties of 
associated Legendre functions we are able to derive various 
expression for the associated Legendre functions with the 
order equal to plus or minus the degree.
\begin{cor}
Let 
$\Re\nu>-\frac12$,
$z\in\mathbb C\setminus(-\infty,1]$.
Then
\begin{eqnarray}
&&P_\nu^\nu(z)=\frac{2^{\nu}\Gamma(\nu+\frac12)}{\sqrt{\pi}}
(z^2-1)^{\frac{\nu}{2}}
+\frac{2^{\nu+1}}{\pi}\sin(\pi\nu)\Gamma(\nu+1)(z^2-1)^{\frac{\nu}{2}}
\int_{z}^\infty
\frac{dw}{(w^2-1)^{\nu+1}}.
\label{Pnunuz}
\end{eqnarray}
\end{cor}
\begin{proof}
Start with the
connection relation
\cite[(14.9.15)]{NIST:DLMF}
\[
P_\nu^\mu(z)
=\frac{\Gamma(\nu+\mu+1)}{\Gamma(\nu-\mu+1)}
P_\nu^{-\mu}(z)
+\frac{2}{\pi}\sin(\pi\mu)\Gamma(\nu+\mu+1)
\biQ_\nu^\mu(z),
\]
then, relying on \eqref{legPnumnu}, the choice 
$\mu=\nu$ completes the proof.
\end{proof}
\begin{rem}
Observe that if $\nu=n\in\mathbb N_0$ 
then
\[
P_n^n(z)=
\frac{2^n\Gamma(n+\frac12)
(z^2-1)^{\frac{n}{2}}}
{\sqrt{\pi}}
=(2n-1)!!(z^2-1)^{\frac{n}{2}},
\]
where we have used
\cite[(6.1.12)]{Abra},
and $(\cdot)!!$ is the
double factorial symbol.
\end{rem}
\begin{cor}
Let 
$\Re\nu>0$, 
$z\in\mathbb C\setminus(-\infty,1]$.
Then
\[
P_{-\nu}^{-\nu}(z)=\frac{1}{2^{\nu-1}\Gamma(\nu)
(z^2-1)^{\frac{\nu}{2}}}
\int_{1}^z (w^2-1)^{\nu-1}{dw}.
\]
\end{cor}
\begin{proof}
{Starting} with \eqref{Qnunuz} and using the Whipple relation for 
associated Legendre functions 
\eqref{Whipple},
followed by the application of the 
involution
\cite[Section 2]{Cohlonpar}
$\zeta(z):=\log\coth\frac{z}{2}$,
then making the change of variables
$w=\zeta/\sqrt{\zeta^2-1}$
completes the proof.
\end{proof}


An interesting definite integral follows
from the behavior of 
the above integral representation near the singularity
at $z=1$.
Using 
\cite[(14.9.15)]{NIST:DLMF},
\yoro{then}
\begin{equation}
P_\nu^{-\nu}(z)=\frac{1}{\Gamma(2\nu+1)}
\left(
P_\nu^\nu(z)-\frac{2 \sin(\pi\nu)}{\pi}
\biQ_\nu^\nu(z)\right).
\end{equation}
After replacement of (\ref{Qnunuz}) and 
(\ref{Pnunuz}) in (\ref{legPnumnu}) we obtain
\begin{eqnarray}\label{2integrals}
&&P_\nu^{-\nu}(z)=-\frac{\sin(\pi\nu)(z^2-1)^{\frac{\nu}{2}}}
{\sqrt{\pi}\,2^{\nu-1}\Gamma(\nu+\frac12)} 
\left(\int_{z/\sqrt{z^2-1}}^\infty \frac{dw}{(w^2-1)^{-\nu+1/2}}
+\int_{z}^\infty \frac{dw}{(w^2-1)^{\nu+1}}\right)
\nonumber\\
&&\hspace{1.2cm}=-\frac{\sin(\pi\nu)(z^2-1)^{\frac{\nu}{2}}}
{\sqrt{\pi}\,2^{\nu-1}\Gamma(\nu+\frac12)} 
\int_1^\infty \frac{dw}{(w^2-1)^{\nu+1}}\nonumber\\
&&\hspace{1.2cm}=\frac{(z^2-1)^{\frac{\nu}{2}}}{2^\nu\Gamma(\nu+1)}.
\end{eqnarray}
Which is simply a re-evaluation of \eqref{legPnumnu}.
From the previous identities the following result follows.

\begin{cor}
\label{Robertocor}
Let $-\frac 12\le \Re \nu\le 0$. Then
\begin{eqnarray*}
&&\frac{\Gamma (-\nu)\Gamma(\nu+\frac{1}{2})}{2\sqrt{\pi}}
=\int_{z/\sqrt{z^2-1}}^\infty \frac{dw}{(w^2-1)^{-\nu+1/2}}
+\int_{z}^\infty \frac{dw}{(w^2-1)^{\nu+1}}
=\int_1^\infty \frac{dw}{(w^2-1)^{\nu+1}}.
\end{eqnarray*}
\end{cor}
\begin{proof}
The formula follows 
after a straightforward calculation starting from \eqref{2integrals}, {making the change of variables
$w=\zeta/\sqrt{\zeta^2-1}$
in the first integral} and
taking into account \eqref{legPnumnu}. 
\end{proof}
\section{Ferrers functions of the first and second kind}

\noindent 
The Ferrers functions
of the first and second
kinds (associated Legendre functions of the first and second \yoro{kinds}
on-the-cut) ${\sf P}_\nu^\mu:(-1,1)\to\mathbb C$ is defined 
in \cite[(14.3.1)]{NIST:DLMF}
\begin{eqnarray}
&&\hspace{-1cm}{\sf P}_\nu^\mu(x):=\left(\frac{1+x}{1-x}\right)^{\frac{\mu}{2}}
\frac{1}{\Gamma(1-\mu)}\Fhyp21{-\nu,\nu+1}{1-\mu}{\frac{1-x}{2}},
\label{FerPdefn}\\
&&\hspace{0.1cm}=\frac{(1-x^2)^{\frac{\mu}{2}}}{2^\mu\Gamma(1+\mu)}
\Fhyp21{\nu+\mu+1,\mu-\nu}{1+\mu}{\frac{1-x}{2}},\label{EulerFerP}
\end{eqnarray}
where we have applied the Euler transformation \eqref{Euler}, to 
the single summation definition
of the Ferrers function of the first kind
produces the second representation for the
Ferrers function of the first kind.
Also, ${\sf Q}_\nu^\mu:(-1,1)\to\mathbb C$ is defined 
in \cite[(14.3.2)]{NIST:DLMF}
\begin{eqnarray}
{\sf Q}_\nu^\mu(x):=\frac{\pi}{2\sin(\pi\mu)}
&\Biggl[&\frac{\cos(\pi\mu)}{\Gamma(1-\mu)}\left(\frac{1+x}{1-x}
\right)^{\frac{\mu}{2}}\Fhyp21{-\nu,\nu+1}{1-\mu}{\frac{1-x}{2}}\nonumber\\
&-&\frac{\Gamma(\nu+\mu+1)}{\Gamma(\nu-\mu+1)}
\left(\frac{1-x}{1+x}\right)^{\frac{\mu}{2}}\frac{1}{\Gamma(1+\mu)}
\Fhyp21{-\nu,\nu+1}{1+\mu}{\frac{1-x}{2}}\Biggr],
\label{FerQdefn}
\end{eqnarray}
where $\mu\not\in\mathbb Z$.
However, ${\sf Q}_\nu^\mu(x)$ can be analytically continued for 
$\mu\in\mathbb Z$ which is demonstrated by \cite[(14.3.12)]{NIST:DLMF},

\begin{eqnarray}
{\sf Q}_\nu^\mu (x)=\frac{\sqrt{\pi}\, 2^{\mu-1}}{(1-x^2)^{\frac{\mu}{2}}} 
&\Bigg[&
\frac{-\sin\left( \tfrac{\pi}{2}(\nu+\mu)\right)
\Gamma(\frac{\nu+\mu+1}{2})}{\Gamma(\frac{\nu-\mu+2}{2})}
\Fhyp21{-\frac{\nu+\mu}2,\frac{\nu-\mu+1}2}{\frac12}{x^2}\nonumber\\
&+&\frac{2\cos\left( \tfrac{\pi}{2}(\nu+\mu)\right)
\Gamma(\frac{\nu+\mu+2}{2})\,x}{\Gamma(\frac{\nu-\mu+1}{2})}
\Fhyp21{\frac{-\nu-\mu+1}{2}, \frac{\nu-\mu+2}{2}}{\frac32}{x^2}\Bigg].
\label{FerQx2}
\end{eqnarray}
Another hypergeometric representation
of the Ferrers function of the second
kind which we will use below is
\begin{eqnarray}
{\sf Q}_\nu^\mu(x)=\frac{2^{\mu-1}\cos(\pi\mu)}{(1-x^2)^{\frac{\mu}{2}}}
&\Gamma(\mu)&x^{\nu+\mu}\Fhyp21{\frac{-\nu-\mu}{2}, \frac{-\nu-\mu+1}{2}}{1-\mu}{\frac{x^2-1}{x^2}}
\nonumber\\ 
&+&\frac{\Gamma(\nu+\mu+1)\Gamma(-\mu)}{2^{\mu+1}
\Gamma(\nu-\mu+1)}(1-x^2)^{\frac{\mu}{2}} x^{\nu-\mu}
\Fhyp21{\frac{\mu-\nu}{2},\frac{\mu-\nu+1}{2}}{\mu+1}{\frac{x^2-1}{x^2}},
\label{repFerQonem1ox2}
\end{eqnarray}
which can be obtained by a limiting procedure \cite[cf.~(14.23.5)]{NIST:DLMF}
\begin{eqnarray*}
\label{FerrersQ}
{\sf Q}_\nu^\mu (x)&=&\frac{\Gamma(\nu+\mu+1)}{2}\left(\expe^{-\frac{1}
{2}i\pi\mu}
\biQ_\nu^\mu (x+i0)
+\expe^{\frac{1}{2}i\pi\mu}
\biQ_\nu^\mu (x-i0)\right)
,
\end{eqnarray*}
for $x\in(-1,1)$, starting from \cite[Entry 29, p.~162]{MOS}.

\subsection{The Ferrers function of the first kind}

Here we derive interesting derivative formulae and integral representations
for the Ferrers function of the first kind.

\medskip
\noindent 
{
First we treat some multi-integrals of the Ferrers function of the first kind from the singularity at $x=1$.
}

\begin{rem}
Using 
\eqref{EulerFerP}
and 
\eqref{derhyp}
produces the following derivative formula for $n\in\mathbb N_0$, 
$x\in(-1,1)$, $\nu,\mu\in\mathbb C$, namely
\begin{equation}
\frac{d^n}{dx^n}
\frac{{\sf P}_\nu^{-\mu}(x)}{(1-x^2)^{\frac{\mu}{2}}}
=(-1)^n (\nu+\mu+1)_n(\mu-\nu)_n
\frac{{\sf P}_\nu^{-\mu-n}(x)}{(1-x^2)^{\frac{\mu+n}{2}}}.
\label{derivnFerPnummu}
\end{equation}
\label{renderFerP}
\end{rem}

\noindent From \eqref{EulerFerP} integral representations 
\yoro{can be obtained}
through repeated integration. For instance, the single integral 
result is given as follows.
\begin{cor}
Let $x\in(-1,1)$,
$\nu,\mu\in\mathbb C$.
Then
\begin{equation}
\int_x^1
\frac{{\sf P}_\nu^{-\mu}(w)}{(1-w^2)^{\frac{\mu}{2}}}
dw
=
\frac{1}{(\nu+\mu)(\nu-\mu+1)}\left(
\frac{1}{2^{\mu-1}\Gamma(\mu)}
-
\frac{{\sf P}_\nu^{-\mu+1}(x)}{(1-x^2)^{\frac{\mu-1}{2}}}\right).
\end{equation}
\end{cor}
\begin{proof}
In order to derive this result,
\yoro{integrate} 
\eqref{derivnFerPnummu}
for $n=1$ with 
the fundamental theorem
of calculus \eqref{funthmcalc}
and taking advantage of
\cite[(14.8.1)]{NIST:DLMF}
\begin{equation}
\label{singFerPlim1}
{\sf P}_\nu^{\mu}(x)\sim\frac{1}{\Gamma(1-\mu)}
\left(\frac{2}{1-x}\right)^{\frac{\mu}{2}},
\end{equation}
as $x\to1^{-}$,
which completes the proof.
\end{proof}
The above result \yoro{can be generalized}
by repeatedly integrating the above
formula.

\begin{thm}
\label{assLegmultint-2}
Let $n\in\mathbb N_0$, $x\in(-1,1)$,
$\nu,\mu\in\mathbb C$.
Then
\begin{eqnarray}
&&\hspace{-2.5cm}\int_x^1\cdots\int_x^1
\frac{{\sf P}_\nu^{-\mu}(w)}{(1-w^2)^{\frac{\mu}{2}}}(dw)^n=
\frac{1}{(-\nu-\mu)_n(\nu-\mu+1)_n}
\Biggl(\frac
{{\sf P}_\nu^{-\mu+n}(x)}
{ (1-x^2)^{\frac{\mu-n}{2}}}\nonumber\\
&&\hspace{0.5cm}-\frac
{(-1)^n(-\mu)_n}
{2^{\mu-n}\Gamma(\mu+1)}
\sum_{k=0}^{n-1}
\frac{(\nu+\mu+1-n)_k(\mu-\nu-n)_k}
{k!(\mu+1-n)_k}
\left(\frac{1-x}{2}\right)^k
\Biggr)
\label{twosumFerP}
\\
&&\hspace{0.5cm}=
\frac{(1-x)^n}{2^\mu n!\Gamma(\mu+1)}
\Fhyp32{\nu+\mu+1,\mu-\nu,1}{\mu+1,n+1}{\frac{1-x}{2}}.
\label{assFerPtwosum}
\end{eqnarray}
\end{thm}
\begin{proof}
Repeated integration of 
\eqref{derivnFerPnummu} while
noting 
\eqref{singFerPlim1},
and using induction 
with with \eqref{Pochdef}
derives the two sum 
expression \eqref{twosumFerP}.
By rewriting the Ferrers function of the first
kind on the right-hand side of \eqref{assFerPtwosum} in terms of the Gauss 
hypergeometric representation
\eqref{EulerFerP}, 
the finite sum term cancels the first $n$ 
terms of the $k$ sum, and rewriting 
the resulting expression
shows it can be written 
in terms of a nonterminating ${}_3F_2$.
This completes the proof.
\end{proof}

\zoro{
\begin{rem}
It is clear that 
Theorem \ref{assLegmultint-2}
is a generalization of \cite[(14.6.6)]{NIST:DLMF}
\begin{equation}
{\sf P}_\nu^{-n}(x)=
(1-x^2)^{-\frac{n}{2}}
\int_x^1\cdots\int_x^1
{\sf P}_\nu(w)(dw)^n,
\end{equation}
by considering the specialization $\mu=0$ in \eqref{twosumFerP}, where we have used 
\cite[(14.9.3)]{NIST:DLMF}
\begin{equation}
{\sf P}_\nu^{-n}(x)=(-1)^n
\frac{\Gamma(\nu-n+1)}
{\Gamma(\nu+n+1)}
{\sf P}_\nu^n(x),\quad n\in\mathbb N_0.
\end{equation}
\end{rem}
}

\medskip
\noindent \yoro{By using} the antiderivative \cite[(14.17.2)]{NIST:DLMF} 
\begin{equation}
\int
\frac{{\sf P}_\nu^{-\mu}(w)}
{(1-w^2)^{\frac{\mu}{2}}}dw
=\frac{1}{(\nu+\mu)(\nu-\mu+1)}
\frac{{\sf P}_\nu^{-\mu+1}(x)}{(1-x^2)^{\frac{\mu-1}{2}}}+C,
\label{intrepFerPDLMF}
\end{equation}
where $C$ is an arbitrary constant 
to derive
some interesting integral representations
for Ferrers functions of the first kind.
Also, 
\yoro{by using} this formula to obtain
a useful derivative formula for Ferrers functions
of the first kind {(see Remark
\ref{renderFerP})}.


\medskip
{\noindent 
Next we treat some multi-integrals of the Ferrers function 
of the first kind from the origin.
}

\medskip 
\noindent 
Evaluation of \eqref{intrepFerPDLMF}
at the endpoints of integration 
produces the following results.

\begin{thm}
\label{intFerPoncezx}
Let $x\in(-1,1)$, $\nu,\mu\in\mathbb C$. Then
\begin{eqnarray}
&&\hspace{-1cm}\int_0^x \frac{{\sf P}_\nu^{-\mu}(w)}
{(1-w^2)^{\frac{\mu}{2}}}dw=
\frac{1}{(\nu+\mu)(\nu-\mu+1)}\left(
\frac{{\sf P}_\nu^{-\mu+1}(x)}{(1-x^2)^{\frac{\mu-1}{2}}}
-\frac{\sqrt{\pi}}{2^{\mu-1}
\Gamma\left(\frac{\nu+\mu+1}{2}\right)
\Gamma\left(\frac{\mu-\nu}{2}\right)
}\right).
\end{eqnarray}
\end{thm}
\begin{proof}
Evaluating \eqref{intrepFerPDLMF}
at the endpoints of integration using 
\cite[(14.5.1)]{NIST:DLMF}
\begin{equation}
{\sf P}_\nu^{-\mu}(0)
=\frac{\sqrt{\pi}}{2^\mu
\Gamma\left(\frac{\nu+\mu+2}{2}\right)
\Gamma\left(\frac{\mu-\nu+1}{2}\right)
},
\label{FerPatzero}
\end{equation}
completes the proof.
\end{proof}

\begin{thm}
Let $n\in\N_0$, $x\in(-1,1)$, $\nu,\mu\in\mathbb C$. Then
\begin{eqnarray}
\int_0^x\cdots\int_0^x \frac{{\sf P}_\nu^{-\mu}(w)}
{(1-w^2)^{\frac{\mu}{2}}}&&(dw)^n=
(-1)^n\sum_{k=n}^\infty (\nu+\mu+1)_{k-n}(\mu-\nu)_{k-n}
{\sf P}_\nu^{-\mu+n-k}(0)\\ =&&
\frac{(-1)^n}{(-\nu-\mu)_n(\nu-\mu+1)_n} \nonumber\\ 
&&
\times\left(\frac{{\sf P}_\nu^{-\mu+n}(x)}{(1-x^2)^{\frac{\mu-n}{2}}}
-\frac{\sqrt{\pi}}{2^{\mu-n}}\sum_{k=0}^{n-1}\frac{(\nu+\mu-n+1)_k
(\mu-\nu-n)_k\left(-\frac{x}{2}\right)^k}
{k!\,\Gamma\left(\frac{\nu+\mu+2-n+k}{2}
\right)\Gamma\left(\frac{\mu-\nu+1-n+k}{2}\right)}\right)\\
=&&
\frac{\sqrt{\pi}\,x^{n}}
{2^\mu n!\Gamma\left(\frac{\nu+\mu+2}{2}\right)
\Gamma\left(\frac{\mu-\nu+1}{2}\right)}
\Fhyp32{\frac{\mu-\nu}{2},\frac{\nu+\mu+1}{2},1}
{\frac{n+1}{2},\frac{n+2}{2}}{x^2}
\nonumber\\
&&-\frac{\sqrt{\pi}\,x^{n+1}}
{2^{\mu-1}(n+1)!\Gamma\left(\frac{\nu+\mu+1}{2}\right)
\Gamma\left(\frac{\mu-\nu}{2}\right)}
\Fhyp32{\frac{\mu-\nu+1}{2},\frac{\nu+\mu+2}{2},1}
{\frac{n+2}{2},\frac{n+3}{2}}{x^2}.
\end{eqnarray}
\label{repintFerPzx}
\end{thm}

\begin{proof}
Repeatedly applying Theorem 
\ref{intFerPoncezx} without evaluating
${\sf P}_\nu^\mu(0)$ 
and then computing the Maclaurin
expansion of ${\sf P}_\nu^{-\mu+n}(x)/(1-x^2)^{(\mu-n)/2}$
yields the first expression.
Using
induction evaluating ${\sf P}_\nu^{-\mu+n}(0)$
with \eqref{Pochdef} 
produces the second expression.
The third expression is obtained by starting
with the first expression, evaluating
${\sf P}_{\nu}^{-\mu+n-k}(0)$, shifting the sum index by $n$ and splitting the sum into even and odd parts.
\end{proof}

\medskip
{
\noindent On the other hand, \yoro{by applying} the antiderivative \cite[(14.17.1)]{NIST:DLMF} 
\begin{equation}
\int
\frac{{\sf P}_\nu^{\mu}(w)}
{(1-w^2)^{\frac{\mu}{2}}}dw
=-
\frac{{\sf P}_\nu^{\mu-1}(x)}{(1-x^2)^{\frac{\mu-1}{2}}}+C,
\label{intrepFerPDLMFb}
\end{equation}
where $C$ is an arbitrary constant 
to derive
some interesting integral representations
for Ferrers functions of the first kind.
Also, \yoro{utilizing} this formula to obtain
a useful derivative formula for Ferrers functions
of the first kind.
}

{
\begin{rem}
Differentiating the above \yoro{result} 
produces the following
formula for $x\in(-1,1)$, $\nu,\mu\in\mathbb C$,
\begin{eqnarray}
\frac{d^n}{dx^n}
\frac{{\sf P}_\nu^{\mu}(x)}{(1-x^2)^{\frac{\mu}{2}}}
=
\frac{(-1)^n{\sf P}_\nu^{\mu+n}(x)}{(1-x^2)^{\frac{\mu+n}{2}}}.
\label{derivFerPnumu}
\end{eqnarray}
\label{rem47}
\end{rem}
}

\medskip
\noindent 
Evaluation of \eqref{intrepFerPDLMFb}
at the endpoints of integration 
produces the following results.

\begin{thm}
Let $x\in(-1,1)$, $\nu,\mu\in\mathbb C$, $\Re\mu>0$. Then
\begin{equation}
\int_x^1 (1-w^2)^{\frac{\mu}{2}} {\sf P}_\nu^{-\mu}(w)dw=
(1-x^2)^{\frac{\mu+1}{2}}{\sf P}_\nu^{-\mu-1}(x).
\end{equation}
\label{repintFerPzxbb}
\end{thm}

\begin{proof}
In order to derive this result
\yoro{integrate}
\eqref{derivFerPnumu}
for $n=1$ with 
the fundamental theorem
of calculus \eqref{funthmcalc}
and taking advantage of
cf.~\eqref{singFerPlim1}
\[
(1-x^2)^{\frac{\mu}{2}}{\sf P}_\nu^{-\mu}(x)\sim 0,
\]
as $x\to1^{-}$.
This completes the proof.
\end{proof}

\begin{thm}
Let $x\in(-1,1)$, $\nu,\mu\in\mathbb C$, $\Re\mu>0$. Then
\begin{eqnarray}
\int_x^1\cdots\int_x^1 
(1-w^2)^{\frac{\mu}{2}}{\sf P}_\nu^{-\mu}(w)
&&(dw)^n=(1-x^2)^{\frac{\mu+n}2}{\sf P}^{-\mu-n}_\nu(x).
\end{eqnarray}
\label{repintFerPzxb}
\end{thm}

\begin{proof}
Repeatedly applying Theorem 
\ref{repintFerPzxbb} through induction
proves the result.
\end{proof}

\begin{thm}
\label{intFerPoncezxb}
Let $n\in\N_0$, $x\in(-1,1)$, $\nu,\mu\in\mathbb C$. Then
\begin{eqnarray}
&&\hspace{-1cm}\int_0^x \frac{{\sf P}_\nu^{\mu}(w)}
{(1-w^2)^{\frac{\mu}{2}}}dw=
-
\frac{{\sf P}_\nu^{\mu-1}(x)}{(1-x^2)^{\frac{\mu-1}{2}}}
+\frac{2^{\mu-1}\sqrt{\pi}}{
\Gamma\left(\frac{\nu-\mu+3}{2}\right)
\Gamma\left(\frac{-\nu-\mu+2}{2}\right)
}.
\end{eqnarray}
\end{thm}
\begin{proof}
Evaluating \eqref{intrepFerPDLMF} at the endpoints of integration using 
\eqref{FerPatzero} completes the proof.
\end{proof}

\begin{thm}
Let $x\in(-1,1)$, $\nu,\mu\in\mathbb C$. Then
\begin{eqnarray}
\int_0^x\cdots\int_0^x \frac{{\sf P}_\nu^{\mu}(w)}
{(1-w^2)^{\frac{\mu}{2}}}&&(dw)^n
=(-1)^n\sum_{k=n}^\infty \frac{(-x)^k}{k!}{\sf P}^{\mu-n+k}_\nu(0)
\nonumber\\ =&&
(-1)^n\left(\frac{{\sf P}_\nu^{\mu-n}(x)}
{(1-x^2)^{\frac{\mu-n}{2}}}
-2^{\mu-n}\sqrt{\pi}\sum_{k=0}^{n-1}\frac{(-2x)^k}
{k!\,\Gamma\left(\frac{\nu-\mu+2+n-k}{2}
\right)\Gamma\left(\frac{-\nu-\mu+1+n-k}{2}\right)}
\right)\nonumber\\
=&&
\frac{\sqrt{\pi}\,2^\mu x^n}{n!\Gamma\left(\frac{\nu-\mu+2}{2}\right)
\Gamma\left(\frac{-\nu-\mu+1}{2}\right)}
\Fhyp32{
\frac{\mu-\nu}{2},\frac{\nu+\mu+1}{2},1}
{\frac{n+1}{2},\frac{n+2}{2}}{x^2}
\nonumber\\
&&
-\frac{\sqrt{\pi}\,2^{\mu+1} x^{n+1}}{(n+1)!\Gamma\left(\frac{\nu-\mu+1}
{2}\right)\Gamma\left(\frac{-\nu-\mu}{2}\right)}\Fhyp32{
\frac{\mu-\nu+1}{2},\frac{\nu+\mu+2}{2},1}
{\frac{n+2}{2},\frac{n+3}{2}}{x^2}.
\end{eqnarray}
\label{repintFerPzxb}
\end{thm}
\begin{proof}
Repeatedly applying Theorem \ref{intFerPoncezxb} without evaluating
${\sf P}_\nu^\mu(0)$ and then computing the Maclaurin
expansion of ${\sf P}_\nu^{\mu-n}(x)/(1-x^2)^{(\mu-n)/2}$
yields the first expression. Using induction evaluating 
${\sf P}_\nu^{\mu-n}(0)$ with \eqref{Pochdef} 
produces the second expression. The third expression is obtained 
by starting with the first expression, evaluating 
${\sf P}_{\nu}^{\mu-n+k}(0)$, shifting the sum index by $n$ and 
splitting the sum into even and odd parts.
\end{proof}


A definite-integral result near the singularity at $x=1$
follows using (\ref{Ferrerssecond}), (\ref{intrepFerrers2numnu}), 
and \eqref{Gauss}, namely
\[
\int_0^1(1-w^2)^{\nu-1}dw=\frac{\sqrt{\pi}\,\Gamma(\nu)}
{2\Gamma(\nu+\tfrac12)},
\]
for $\Re\nu>0$.
The well-known special value (see \cite[(14.5.18)]{NIST:DLMF})
\begin{equation}
{\sf P}_\nu^{-\nu}(x)=\frac{(1-x^2)^{\frac{\nu}{2}}}{2^\nu\Gamma(\nu+1)},
\label{Pnumnu}
\end{equation}
in conjunction with \cite[(8.737.1)]{Grad}, 
yields the following integral {representation}.
\begin{cor} Let $\nu\in\mathbb C$, $x\in(-1,1)$. Then
\begin{equation}
{\sf P}_\nu^\nu(x)=\frac{2^\nu(1-x^2)^{\frac{\nu}{2}}}{\sqrt{\pi}}
\left(
\Gamma(\nu+\tfrac12)
\cos(\pi\nu)
+
\frac{2\Gamma(\nu+1)}{\sqrt{\pi}}
\sin(\pi\nu)
\int_0^x\frac{dw}{(1-w^2)^{\nu+1}}
\right).
\end{equation}
\end{cor}
\begin{proof}
Start with \cite[(14.9.2)]{NIST:DLMF}
\[
{\sf P}_\nu^\mu(x)=\cos(\pi\mu)
\frac{\Gamma(\nu+\mu+1)}{\Gamma(\nu-\mu+1)}
{\sf P}_\nu^{-\mu}(x)
+
\frac{2}{\pi}\sin(\pi\mu)
\frac{\Gamma(\nu+\mu+1)}{\Gamma(\nu-\mu+1)}
{\sf Q}_\nu^{-\mu}(x),
\]
replace $\mu=\nu$, 
then using \eqref{intrepFerrers2numnu}, \eqref{Pnumnu} 
completes the proof.
\end{proof}

\begin{rem}
Note that if $\nu=n\in\mathbb N_0$ 
then
\begin{equation}
{\sf P}_n^n(x)=
\frac{(-2)^n\Gamma(n+\frac12)
(1-x^2)^{\frac{n}{2}}}
{\sqrt{\pi}}
=(-1)^n(2n-1)!!(1-x^2)^{\frac{n}{2}},
\end{equation}
where we have used
\cite[(6.1.12)]{Abra}.
\end{rem}

\subsection{The Ferrers function of the
second kind}

The Ferrers function of the second kind (associated Legendre 
function of the second kind on-the-cut) 
${\sf Q}_\nu^\mu:(-1,1)\to\mathbb C$ is defined
in \eqref{FerQdefn}.

\medskip
\noindent First we treat some multi-integrals of the Ferrers 
function of the second kind to the singularity at $x=1$.

\begin{lem}
Let $x\in(-1,1)$,
$\nu,\mu\in\mathbb C$ 
such that $\mu\not\in-\N$, $\nu-\mu\not\in-\mathbb N_0$.
Then
\begin{eqnarray}
\int_x^1 (1-w^2)^{\frac{\mu}{2}}
{\sf Q}_\nu^{-\mu}(w)&&dw
=(1-x^2)^{\frac{\mu+1}{2}}
{\sf Q}_\nu^{-\mu-1}(x)
-\frac{2^\mu \Gamma(\mu+1)\Gamma(\nu-\mu)}
{\Gamma(\nu+\mu+2)}.
\end{eqnarray}
\label{reflemQnummu}
\end{lem}
\begin{proof}
The Ferrers function
of the second kind as $x$ approaches
the singularity at $x=1$ has the following
behavior
\cite[(14.8.6)]{NIST:DLMF}
\begin{equation}
(1-x^2)^{\frac{\mu}{2}} {\sf Q}_\nu^{-\mu}(x)
\sim
\frac{2^{\mu-1}\Gamma(\mu)\Gamma(\nu-\mu+1)}{\Gamma(\nu+\mu+1)},
\label{FerQapproxonea}
\end{equation}
as $x\to 1^{-}$, { $\Re\mu>0$.
Evaluating \cite[(14.17.1)]{NIST:DLMF} using 
the Ferrers function of the second kind at the} {endpoints of integration noting the above behavior at $x\approx 1$ completes the proof.}
\end{proof}

\begin{rem}
\label{rem415}
Applying the fundamental theorem of calculus
\eqref{funthmcalc}
to Lemma \ref{reflemQnummu}
produces the following derivative formula
for $x\in(-1,1)$,
$\nu,\mu\in\mathbb C$, namely
\begin{equation}
\frac{d^n}{dx^n} (1-x^2)^{\frac{\mu}{2}}
{\sf Q}_\nu^{-\mu}(x)
=(-1)^n (1-x^2)^{\frac{\mu-n}{2}}
{\sf Q}_\nu^{-\mu+n}(x).
\end{equation}
\end{rem}

{
\begin{thm}Let $n\in\N_0$, $x\in(-1,1)$,
$\nu,\mu\in\mathbb C$.
Then
\begin{eqnarray}
&&\hspace{-0.7cm}\int_x^1\cdots\int_x^1
(1-w^2)^{\frac{\mu}{2}}
{\sf Q}_\nu^{-\mu}(w)
(dw)^n=
(1-x^2)^{\frac{\mu+n}{2}}
{\sf Q}_{\nu}^{-\mu-n}(x)\nonumber\\
&&
\label{cor3.14a}
\hspace{-0.4cm}-
\frac{(-1)^n 2^{\mu+n-1}\Gamma(\mu)
\Gamma(\nu-\mu+1)(\mu)_n}
{\Gamma(\nu+\mu+1)(\mu-\nu)_n(\nu+\mu+1)_n}
\sum_{k=0}^{n-1}
\frac
{(\nu-\mu-n+1)_k(-\nu-\mu-n)_k}
{k!(-\mu-n+1)_k}
\left(
\frac{1-x}{2}
\right)^k\\
&&=-\frac{\pi}{2}\cot(\pi\mu)(1-x^2)^{\frac{\mu+n}{2}}{\sf P}_\nu^{-\mu-n}(x)\nonumber\\
&&\hspace{1cm}
+\frac{2^{\mu-1}\Gamma(\mu)\Gamma(\nu-\mu+1)(1-x)^n}
{n!\Gamma(\nu+\mu+1)}
\Fhyp32{\nu-\mu+1,-\nu-\mu,1}{n+1,1-\mu}{\frac{1-x}{2}}.
\end{eqnarray}
\label{CorFerQnummu}
\end{thm}
\begin{proof}
Repeatedly applying Lemma \ref{reflemQnummu} to itself
using induction
with \eqref{Pochdef} produces the first formula.
The second formula is obtained by rewriting the finite sum 
as a sum from $0$ to $\infty$ and subtracting the sum from 
$n$ to $\infty$, and then finally utilizing \eqref{FerQdefn}.
\end{proof}
}

\begin{rem}
{Taking the $\mu\to 0$ limit} in Corollary 
\ref{CorFerQnummu}
produces the following multi-integration
result for $x\in(-1,1)$, $\nu,\mu\in\mathbb C$, namely
\begin{eqnarray}
&&
\int_x^1\cdots\int_x^1
{\sf Q}_\nu(w)(dw)^n=(1-x^2)^{\frac{n}{2}}
{\sf Q}_\nu^{-n}(x)\nonumber\\
&&\hspace{4.0cm}
-\frac{(-1)^n2^{n-1}(n-1)!}{(-\nu)_n
(\nu+1)_n}
\sum_{k=0}^{n-1}
\frac{(\nu-n+1)_k(-\nu-n)_k}
{k! (1-n)_k}
\left(\frac{1-x}{2}\right)^k.
\end{eqnarray}
\label{FerQlimzero1}
\end{rem}

\noindent 
{
Now we present a similar result for the Ferrers function
of the second kind with order $\mu$ instead of $-\mu$.
\begin{lem}
Let  $x\in(-1,1)$,
$\nu,\mu\in\mathbb C$ 
such that $\mu\not\in-\N$, $\nu-\mu\not\in-\mathbb N_0$.
Then
\begin{equation}
\int_x^1 (1-w^2)^{\frac{\mu}{2}}
{\sf Q}_\nu^{\mu}(w)dw
=\frac{1}{(\mu-\nu) (\nu+\mu+1)}\left((1-x^2)^{\frac{\mu+1}{2}}
{\sf Q}_\nu^{\mu+1}(x)
+\frac{2^\mu\pi \Gamma(\mu+1)}
{\Gamma(-\mu-\frac32)\Gamma(\mu+\frac52)}\right).
\end{equation}
\label{reflemQnumu}
\end{lem}
\begin{proof}
The Ferrers function
of the second kind as $x$ approaches
the singularity at $x=1$ has the following
behavior
(cf.~\cite[(14.8.4)]{NIST:DLMF})
\begin{equation}
(1-x^2)^{\frac{\mu}{2}} {\sf Q}_\nu^{\mu}(x)
\sim
-\frac{2^{\mu-1}\pi\Gamma(\mu)}{\Gamma(-\mu-\frac12)\Gamma(\mu+\frac32)},
\label{FerQapproxoneb}
\end{equation}
as $x\to 1^{-}$, $\Re\mu>0$.
Evaluating \cite[(14.17.2)]{NIST:DLMF} using 
the Ferrers function of the second kind at the endpoints of integration
noting the above behavior at $x\approx 1$ completes the proof.
\end{proof}
}

\begin{rem}
Applying the fundamental theorem of calculus
\label{rem419}
\eqref{funthmcalc}
to Lemma \ref{reflemQnumu}
produces the following derivative formula
for $x\in(-1,1)$,
$\nu,\mu\in\mathbb C$, namely
\begin{equation}
\frac{d^n}{dx^n} (1-x^2)^{\frac{\mu}{2}}
{\sf Q}_\nu^{\mu}(x)
=(-1)^n (\nu-\mu+1)_n(-\nu-\mu)_n (1-x^2)^{\frac{\mu-n}{2}}
{\sf Q}_\nu^{\mu-n}(x).
\end{equation}
\end{rem}

\begin{thm}Let $n\in\N_0$, $x\in(-1,1)$,
$\nu,\mu\in\mathbb C$.
Then
\begin{eqnarray}
&&\hspace{-0.7cm}\int_x^1\cdots\int_x^1
(1-w^2)^{\frac{\mu}{2}}
{\sf Q}_\nu^{\mu}(w)
(dw)^n=\frac{1}{(\mu-\nu)_n(\nu+\mu+1)_n}
\Biggl(
(1-x^2)^{\frac{\mu+n}{2}}
{\sf Q}_{\nu}^{\mu+n}(x)\nonumber\\
&&
\label{cor3.14b}
\hspace{-0.4cm}+
(-1)^{n-1} 2^{\mu+n-1}\cos(\pi\mu)\Gamma(\mu+n)
\sum_{k=0}^{n-1}
\frac
{(\nu-\mu-n+1)_k(-\nu-\mu-n)_k}
{k!(-\mu-n+1)_k}
\left(
\frac{1-x}{2}
\right)^k
\Biggr)\\
&&=-\frac{\pi\Gamma(\nu+\mu+1)}{2\Gamma(\nu-\mu+1)\sin(\pi\mu)}
(1-x^2)^{\frac{\mu+n}{2}}{\sf P}_\nu^{-\mu-n}(x)\nonumber\\
&&\hspace{1cm}
+\frac{2^{\mu-1}\cos(\pi\mu)\Gamma(\mu)}
{n!}(1-x)^n
\Fhyp32{\nu-\mu+1,-\nu-\mu,1}{n+1,1-\mu}{\frac{1-x}{2}}.
\end{eqnarray}
\label{CorFerQnumu}
\end{thm}
\begin{proof}
Repeatedly applying Lemma \ref{reflemQnumu} to itself
using induction
with \eqref{Pochdef} 
produces the first formula.
The second formula is obtained by rewriting
the finite sum as a sum from $0$ to $\infty$ and subtracting the sum 
from $n$ to $\infty$, and then finally utilizing 
\eqref{FerQdefn}.
\end{proof}

{
\begin{rem}
\yoro{An} interesting discussion is concerning whether  
Lemma \ref{Robertolem} \yoro{might be used} to obtain new generalized 
hypergeometric representations for Corollaries \ref{CorFerQnummu},
\ref{CorFerQnumu}. In order to do this, one must compute the 
one-sided Taylor
expansions of the relevant functions about the singular point $x=1$
(the relevant functions are well-behaved at this singular point).
This is readily possible, but is not practical due to the fact
that the behavior of the functions in question near the
singularity changes in form depending on whether $\Re\mu{\lessgtr} 0$
(see \eqref{FerQapproxonea}, \eqref{FerQapproxoneb}). 
The derivative terms in the Taylor series necessarily cross the 
$\Re\mu=0$ boundary, so a simple result from this Lemma does not 
seem to be practical.
\end{rem}
}

{
\begin{rem}
Taking $\mu\to0$ limit in Corollary \ref{CorFerQnumu}
produces the following multi-integration result for 
$x\in(-1,1)$, $\nu,\mu\in\mathbb C$, namely
\begin{eqnarray}
&&
\int_x^1\cdots\int_x^1
{\sf Q}_\nu(w)(dw)^n=
\frac{1}{(-\nu)_n(\nu+1)_n}
\Biggl((1-x^2)^{\frac{n}{2}}
{\sf Q}_\nu^{n}(x)\nonumber\\
&&\hspace{4.0cm}
+(-1)^{n-1}2^{n-1}(n-1)!
\sum_{k=0}^{n-1}
\frac{(\nu-n+1)_k(-\nu-n)_k}
{k! (1-n)_k}
\left(\frac{1-x}{2}\right)^k\Biggr).
\end{eqnarray}
\label{FerQlimzero2}
\end{rem}
}

{
\begin{rem}
Note that in Corollaries \ref{CorFerQnummu}, 
\ref{CorFerQnumu}, it is tempting to consider the 
$\mu\to0$ limit using their ${}_3F_2$ representations.
However, to zeroth order in $\mu$, the limits cancel. 
One must then determine a first order approximation in 
$\mu$ to determine the limit behavior. 
\yoro{After performing} 
this calculation in both of these
situations, then it turns out that \yoro{the} result is given in terms
of the sum of several double hypergeometric series of 
Kamp\'{e} de F\'{e}rier type
\yoro{(see e.g., \cite[p.~27]{SriKarl})}. Since this result is very cumbersome
and doesn't really shed much light on these limits, we have instead
presented the above Remarks \ref{FerQlimzero1},
\ref{FerQlimzero2}. It should also be pointed out
regarding the fact that there does not seem to be an analogous 
formula for the Ferrers function of the second kind in the 
classical list \cite[(14.6.6-8)]{NIST:DLMF}, the above reasoning
most likely explains this fact. The conclusion is that 
any formula for negative integer order Ferrers functions of 
the second kind will involve a finite sum of polynomial terms
in addition to the multi-integral of ${\sf Q}_\nu(x)$, as indicated
in Remarks \ref{FerQlimzero1},
\ref{FerQlimzero2}.
\end{rem}
}

\medskip
\noindent
{
Next we treat some multi-integrals of the 
Ferrers function of the second kind from the origin.
}

{
{
\begin{rem}
\label{rem424}
Applying the fundamental theorem of calculus
\eqref{funthmcalc}
to \cite[(14.17.2)]{NIST:DLMF}
produces the following derivative formula
for $n\in\N_0$, $x\in(-1,1)$,
$\nu,\mu\in\mathbb C$, namely
\begin{equation}
\frac{d^n}{dx^n} 
\frac{{\sf Q}_\nu^{-\mu}(x)}
{(1-x^2)^{\frac{\mu}{2}}}
=(-1)^n (\nu+\mu+1)_n(\mu-\nu)_n 
\frac{{\sf Q}_\nu^{-\mu-n}(x)}
{(1-x^2)^{\frac{\mu+n}{2}}}.
\end{equation}
\end{rem}
}
}

\begin{thm}
Let $x\in(-1,1)$, $\nu,\mu\in\mathbb C$.
Then
\begin{eqnarray}
&&\hspace{-1cm}\int_0^x\frac
{{\sf Q}_\nu^{-\mu}(w)}
{(1-w^2)^{\frac{\mu}{2}}}
dw
=
\frac{1}{(\nu+\mu)(\nu-\mu+1)}\left(\frac{{\sf Q}_\nu^{-\mu+1}(x)}
{(1-x^2)^{\frac{\mu-1}{2}}}
-\frac{\pi^\frac32
\Gamma\left(\frac{\nu-\mu+2}{2}\right)
}
{2^{\mu}
\Gamma\left(\frac{\nu+\mu+1}{2}\right)
\Gamma\left(\frac{\mu-\nu-1}{2}\right)
\Gamma\left(\frac{\nu-\mu+3}{2}\right)
}\right).
\end{eqnarray}
\label{thmintQnummu}
\end{thm}
\begin{proof}
Evaluating \cite[{(14.17.2)}]{NIST:DLMF}
(expressed as a Ferrers function
of the second kind)
at the endpoints of integration using 
\cite[(14.5.3)]{NIST:DLMF}
\begin{equation}
\yoro{
{\sf Q}_\nu^{-\mu}(0)
=\frac{\pi^{\frac32}
\Gamma\left(\frac{\nu-\mu+1}{2}\right)
}{2^{\mu+1}
\Gamma\left(\frac{\nu+\mu+2}{2}\right)
\Gamma\left(\frac{\mu-\nu}{2}\right)
\Gamma\left(\frac{\nu-\mu+2}{2}\right)
}},
\label{FerQatzero}
\end{equation}
completes the proof.
\end{proof}

{
\begin{thm}
Let $n\in\N_0$, $x\in(-1,1)$, $\nu,\mu\in\mathbb C$.
Then
\begin{eqnarray}
&&\hspace{-1.3cm}\int_0^x\cdots\int_0^x
\frac{{\sf Q}_\nu^{-\mu}(w)}
{(1-w^2)^{\frac{\mu}{2}}}
(dw)^n
{
=
(-1)^n
\sum_{k=n}^\infty \frac{(-x)^k}{k!}
(\nu+\mu+1)_{k-n}(\mu-\nu)_{k-n}
{\sf Q}_\nu^{-\mu+n-k}(0)}\\
&&=\frac{(-1)^n}{(-\nu-\mu)_n(\nu-\mu+1)_n}\left(
\frac{{\sf Q}_\nu^{-\mu+n}(x)}
{(1-x^2)^{\frac{\mu-n}{2}}}\right.
\nonumber
\\
&&\hspace{1cm}-\left.\frac{\pi^\frac32}
{2^{\mu-n+1}}
\sum_{k=0}^{n-1}
\frac{(\nu+\mu+1-n)_k
(\mu-\nu-n)_k
\Gamma\left(
\frac{\nu-\mu+1+n-k}{2}
\right)\left(-\frac{x}{2}\right)^k
}
{k!\,\Gamma\left(
\frac{\nu+\mu+2-n+k}{2}
\right)
\Gamma\left(
\frac{\mu-\nu-n+k}{2}
\right)
\Gamma\left(
\frac{\nu-\mu+2+n-k}{2}
\right)}\right)\\
&&{=\frac{\pi^\frac32x^n
\Gamma\left(\frac{\nu-\mu+1}{2}\right)
}
{n!2^{\mu+1}
\Gamma\left(\frac{\nu+\mu+2}{2}\right)
\Gamma\left(\frac{\mu-\nu}{2}\right)
\Gamma\left(\frac{\nu-\mu+2}{2}\right)
}\Fhyp32
{\frac{\nu+\mu+1}{2},
\frac{\mu-\nu}{2},1
}
{\frac{n+1}{2},\frac{n+2}{2}}
{x^2}}\nonumber\\
&&\hspace{0.3cm}{+
\frac{\pi^\frac32x^{n+1}
\Gamma\left(\frac{\nu-\mu+2}{2}\right)
}
{(n+1)!2^{\mu}
\Gamma\left(\frac{\nu+\mu+1}{2}\right)
\Gamma\left(\frac{\mu-\nu+1}{2}\right)
\Gamma\left(\frac{\nu-\mu+1}{2}\right)
}
\Fhyp32{\frac{\nu+\mu+2}{2},\frac{\mu-\nu+1}{2},1}
{\frac{n+2}{2},\frac{n+3}{2}}{x^2}}.
\end{eqnarray}
\label{repintFerQzx}
\end{thm}
\begin{proof}
Repeatedly applying Theorem 
\ref{thmintQnummu} without evaluating
${\sf Q}_\nu^{-\mu}(0)$ 
and then computing the Maclaurin
expansion of ${\sf Q}_\nu^{-\mu+n}(x)/(1-x^2)^{(\mu-n)/2}$
yields the first expression.
Using
induction evaluating ${\sf Q}_\nu^{-\mu+n}(0)$
with \eqref{Pochdef} produces the second expression.
The third expression is obtained by starting
with the first expression, evaluating
${\sf Q}_{\nu}^{-\mu+n-k}(0)$, shifting the sum index by $n$ 
and splitting the sum into even and odd parts.
\end{proof}
}

{
\begin{rem}
{Taking the limit as $\mu=0$} in Theorem
\ref{repintFerQzx} 
produces the following multi-integration result,
namely for $x\in(-1,1)$, $\nu,\mu\in\mathbb C$, 
then
\begin{eqnarray}
&&\hspace{-0.5cm}\int_0^x\cdots\int_0^x
{\sf Q}_\nu(w)
(dw)^n
=\frac{1}{(-\nu)_n(\nu+1)_n}
\Biggl((-1)^n
(1-x^2)^{\frac{n}{2}}
{\sf Q}_\nu^{n}(x)
\nonumber
\\
&&\hspace{3cm}+(-1)^{n+1}2^{n-1}\pi^\frac32
\sum_{k=0}^{n-1}
\frac{
(\nu+1-n)_k
(-\nu-n)_k
\Gamma\left(
\frac{\nu+1+n-k}{2}
\right)\left(-\frac{x}{2}\right)^k
}
{k!\,
\Gamma\left(
\frac{\nu+2-n+k}{2}
\right)
\Gamma\left(
\frac{-\nu-n+k}{2}
\right)
\Gamma\left(
\frac{\nu+2+n-k}{2}
\right)
}\Biggr).
\end{eqnarray}
\label{repcorintFerQzx}
\end{rem}
}


\medskip
{
\noindent On the other hand, \yoro{by applying} the antiderivative 
\cite[(14.17.1)]{NIST:DLMF} 
\begin{equation}
\int
\frac{{\sf Q}_\nu^{\mu}(w)}
{(1-w^2)^{\frac{\mu}{2}}}dw
=-
\frac{{\sf Q}_\nu^{\mu-1}(x)}{(1-x^2)^{\frac{\mu-1}{2}}}+C,
\label{intrepFerQDLMFb}
\end{equation}
where $C$ is an arbitrary constant 
to derive
some interesting integral representations for Ferrers functions of the second kind. 
}

\medskip
\noindent {
An examination of the above formula
\eqref{intrepFerQDLMFb} produces the following results.  
For instance, 
\yoro{by applying} this formula to obtain
a useful derivative formula for Ferrers functions
of the second kind.
}

{
\begin{rem}
\label{rem428}
Differentiating the above formula with \eqref{Pochdef} produces 
the following formula for $x\in(-1,1)$, $\nu,\mu\in\mathbb C$,
\begin{eqnarray}
\frac{d^n}{dx^n}\frac{{\sf Q}_\nu^{\mu}(x)}{(1-x^2)^{\frac{\mu}{2}}}
=\frac{(-1)^n{\sf Q}_\nu^{\mu+n}(x)}{(1-x^2)^{\frac{\mu+n}{2}}}.
\end{eqnarray}
\end{rem}
}

{
\begin{thm}
\label{thmeqintFerQoncezxb}
Let $x\in(-1,1)$, $\nu,\mu\in\mathbb C$, such that
$\nu+\mu\not\in-2\mathbb N_0$. Then
\begin{eqnarray}
&&\hspace{-1cm}\int_0^x \frac{{\sf Q}_\nu^{\mu}(w)}
{(1-w^2)^{\frac{\mu}{2}}}dw=
-
\frac{{\sf Q}_\nu^{\mu-1}(x)}{(1-x^2)^{\frac{\mu-1}{2}}}
-\frac{2^{\mu-2}\pi^\frac32\Gamma\left(\frac{\nu+\mu}{2}\right)}{
\Gamma\left(\frac{\nu-\mu+3}{2}\right)
\Gamma\left(\frac{-\nu-\mu-1}{2}\right)
\Gamma\left(\frac{\nu+\mu+3}{2}\right)
}.
\label{eqintFerQoncezxb}
\end{eqnarray}
\end{thm}
\begin{proof}
Evaluating \eqref{intrepFerQDLMFb}
at the endpoints of integration using 
\eqref{FerQatzero}
completes the proof.
\end{proof}
}

{
\begin{thm}
Let $n\in\N_0$, $x\in(-1,1)$, $\nu,\mu\in\mathbb C$. Then
\begin{eqnarray}
&&\hspace{-0.8cm}\int_0^x\cdots\int_0^x \frac{{\sf Q}_\nu^{\mu}(w)}
{(1-w^2)^{\frac{\mu}{2}}}(dw)^n
=(-1)^n\sum_{k=n}^\infty\frac{(-x)^k}{k!}{\sf Q}_\nu^{\mu-n+k}(0)\\
&&=
\frac{(-1)^n {\sf Q}_\nu^{\mu-n}(x)}
{(1-x^2)^{\frac{\mu-n}{2}}}
+\frac{(-1)^n\pi^\frac32}{2^{n+1-\mu}}
\sum_{k=0}^{n-1}\frac{
\Gamma\left(\frac{\nu+\mu+1-n+k}{2}
\right)(-2x)^k
}{k!\,\Gamma\left(\frac{\nu-\mu+2+n-k}{2}
\right)\Gamma\left(\frac{-\nu-\mu-2+n-k}{2}\right)
\Gamma\left(\frac{\nu+\mu+4-n+k}{2}\right)}
\nonumber\\
&&{
\hspace{0.2cm}=\frac{\pi^\frac32\,2^{\mu-1} x^n\Gamma\left(\frac{\nu+\mu+1}{2}\right)
}{
n!\Gamma\left(\frac{\nu-\mu+2}{2}\right)
\Gamma\left(\frac{-\nu-\mu}{2}\right)
\Gamma\left(\frac{\nu+\mu+2}{2}\right)
}
\Fhyp32{
\frac{\mu-\nu}{2},\frac{\nu+\mu+1}{2},1}
{\frac{n+1}{2},\frac{n+2}{2}}{x^2}}
\nonumber\\
&&\hspace{1.5cm}{-\frac{\pi^\frac32\,2^{\mu} x^{n+1}\Gamma\left(\frac{\nu+\mu+2}{2}\right)}{
(n+1)!\Gamma\left(\frac{\nu-\mu+1}{2}\right)
\Gamma\left(\frac{-\nu-\mu-1}{2}\right)
\Gamma\left(\frac{\nu+\mu+3}{2}\right)
}
\Fhyp32{
\frac{\mu-\nu+1}{2},\frac{\nu+\mu+2}{2},1}
{\frac{n+2}{2},\frac{n+3}{2}}{x^2}
}.
\end{eqnarray}
\label{repintFerQzxb}
\end{thm}
}

{
\begin{proof}
Repeatedly applying Theorem 
\ref{thmeqintFerQoncezxb} without evaluating
${\sf Q}_\nu^\mu(0)$ 
and then computing the Maclaurin
expansion of ${\sf Q}_\nu^{\mu-n}(x)/(1-x^2)^{(\mu-n)/2}$
yields the first expression.
Using
induction evaluating ${\sf Q}_\nu^{\mu-n}(0)$
with \eqref{Pochdef} 
produces the second expression.
The third expression is obtained by starting
with the first expression, evaluating
${\sf Q}_{\nu}^{\mu-n+k}(0)$, shifting the sum index by $n$ and 
splitting the sum into even and odd parts.
\end{proof}
}

\begin{thm}
Let 
$\nu\in\mathbb C$.
Then, 
\begin{equation}
\int\frac{dx}{(1-x^2)^{\nu+1}}=x
\,{}_2F_1\left(\begin{array}{c}
\frac12,\nu+1\\[0.2cm]
\frac32\end{array};{x^2}\right)+C
=\frac{2^\nu\Gamma(\nu+\tfrac12)}{\sqrt{\pi}\,(1-x^2)^{\frac{\nu}{2}}}
{\sf Q}_\nu^{-\nu}(x)+C,
\end{equation}
{where $C$ is an arbitrary constant.}
\label{antiderivthmsecond}
\end{thm}

{
\begin{proof}
The Gauss hypergeometric function in the
antiderivative follows using 
(\ref{derhyp}),
(\ref{sechypcontig}),
(\ref{thirdhypcontig}),
as in the proof of Theorem \ref{thmantifirst},
with the Ferrers function following 
directly using
\begin{equation}
{\sf Q}_\nu^{-\nu}(x)=\frac{\sqrt{\pi}\,x(1-x^2)^{\frac{\nu}{2}}}
{2^\nu\Gamma(\nu+\tfrac12)}\,
{}_2F_1\left(
\begin{array}{c}
\frac12,\nu+1\\[0.2cm]
\frac32
\end{array};x^2
\right),
\label{Ferrerssecond}
\end{equation}
which follows from 
\eqref{FerQx2}. 
This completes the 
proof.
\end{proof}
}

\medskip
The following {very simple} 
integral representation for the Ferrers function of the second kind
is a consequence of Theorem \ref{antiderivthmsecond}.
\begin{cor} Let $\nu\in\mathbb C$, $x\in(-1,1)$. Then
\begin{equation}
{\sf Q}_\nu^{-\nu}(x)=\frac{\sqrt{\pi}(1-x^2)^{\frac{\nu}{2}}}
{2^\nu\,\Gamma(\nu+\frac12)}
\int_{0}^x\frac{dw}{(1-w^2)^{\nu+1}}.
\label{intrepFerrers2numnu}
\end{equation}
\label{corQnumnuint}
\end{cor}
\begin{proof}
Evaluating the antiderivative Theorem 
\ref{antiderivthmsecond} at the endpoints of integration 
completes the proof.
\end{proof}

\begin{rem}
One can also show that 
Corollary \ref{corQnumnuint}
also follows directly from Theorem 
\ref{thmintQnummu}.
This is true even though 
Theorem \ref{thmintQnummu}
is not strictly valid for 
$\mu=-\nu$.
The result can be obtained by
taking the limit as $\mu\to-\nu$
in Theorem 
\ref{thmintQnummu}
and using the Gauss hypergeometric
representation of the Ferrers
function of the second kind 
with argument $(1-x^{-2})$, namely
\eqref{repFerQonem1ox2}.
\end{rem}

\begin{cor} Let $\nu\in\mathbb C$, $x\in(-1,1)$. Then
\begin{eqnarray}
&&
{\sf Q}_\nu^\nu(x)=-2^{\nu-1}\sqrt{\pi}\,\Gamma(\nu+\tfrac12)
\sin(\pi\nu)(1-x^2)^{\frac{\nu}{2}}
\nonumber\\[0.2cm]&&\hspace{5.0cm}
+2^\nu\Gamma(\nu+1)\cos(\pi\nu)(1-x^2)^{\frac{\nu}{2}}
\int_0^x\frac{dw}{(1-w^2)^{\nu+1}}.
\end{eqnarray}
\label{corlast}
\end{cor}
\begin{proof}
Using the connection relation \cite[p.~170]{MOS}
\[
{\sf Q}_\nu^{-\mu}(x)=
\frac{\Gamma(\nu-\mu+1)}{\Gamma(\nu+\mu+1)}
\left(
\cos(\pi\mu){\sf Q}_\nu^\mu(x)+\frac{\pi}{2}
\sin(\pi\mu)
{\sf P}_\nu^\mu(x)
\right),
\]
setting $\mu=\nu$, and using the above results 
completes the proof.
\end{proof}

\begin{rem}
Note that if $\nu=n+\frac12$, $n\in\mathbb N_0$ 
then Corollary \ref{corlast}
reduces to the following special value
\begin{equation}
{\sf Q}_{n+\frac12}^{n+\frac12}(x)=
(-1)^{n+1} 2^{n-\frac12} n! \sqrt{\pi}\, (1-x^2)^{\frac{n}2+\frac14}.
\end{equation}
\end{rem}

\section{\yoro{Conclusion}}

\zoro{
In this paper, we explore some
implications of the existence of multi-derivative formulae for associated Legendre functions of the first and second kinds $P_\nu^\mu$, ${\bf Q}_\nu^\mu$, and Ferrers functions of the first and second kinds ${\sf P}_\nu^\mu$, ${\sf Q}_\nu^\mu$.
These multi-derivative formulae (see Remarks \ref{rem39}, \ref{rem315}, \ref{renderFerP}, \ref{rem47}, \ref{rem415},
\ref{rem419}, \ref{rem424}, \ref{rem428}) have the useful property that the degree $(\nu)$ is left unchanged 
by the multi-derivative.
The order
$(\mu)$ is then shifted by unit increments depending on the number of derivatives. These
multi-derivative formulae generalize some classical  multi-derivative formulae for these functions with integer order
\cite[(14.6.1)--(14.6.5)]{NIST:DLMF}.
}

\zoro{
Due to the existence of these multi-derivative formulae, and certain special known values and limiting behaviors near the singularities of these functions, we derive several multi-integral
representations for these functions. 
These multi-integral representations
are shown to be given either in terms of ({\sf i}) a sum of two ${}_3F_2$'s (Theorems \ref{repintFerPzx}, \ref{repintFerPzxb}, 
\ref{repintFerQzx}, \ref{repintFerQzxb}); ({\sf ii}) a ${}_2F_1$ and a ${}_3F_2$ (Theorems \ref{CorFerQnummu},
\ref{CorFerQnumu}); ({\sf iii}) 
a single ${}_3F_2$ 
(Theorems \ref{multintLegP}, \ref{assLegmultint-2}); or ({\sf iv}) a single 
${}_2F_1$ (Theorems \ref{critthmOlverQ},
\ref{repintFerPzxb}).
 These
multi-integral representations generalize some classical multi-integrals for these functions with integer order
\cite[(14.6.6)--(14.6.8)]{NIST:DLMF}.
}

\zoro{
As mentioned in the introduction, many of the functions encountered in this work represent fundamental solutions for the Laplace-Beltrami operator on Riemannian manifolds of constant curvature. Multi-integrals and derivatives of these functions are essential in performing global analysis for these fundamental
solutions on these manifolds. One interesting open problem where this work is almost certainly essential is for obtaining fundamental solutions of natural powers of the Laplace-Beltrami
operator (polyharmonic) on these manifolds. This analysis will be investigated in future publications.}

\vspace{6pt} 
\authorcontributions{H.S.C. and R.S.C.-S. conceived the mathematics;
H.S.C. and R.S.C.-S. wrote the paper.}
\funding{{The research of R.S.C.-S. was funded by Agencia Estatal de 
Investigaci\'on of Spain, grant number PGC-2018-096504-B-C33.}}


\conflictsofinterest{The authors declare no conflict of interest.}
\reftitle{References}


\end{document}